\documentclass[twoside,english,links,13pt,a4paper]{article} 
\pdfoutput=1
\usepackage[lmargin=3cm,tmargin=3cm,bmargin=3cm,rmargin=3cm]{geometry}
\usepackage[numbers,square,sort,comma]{natbib}
\usepackage{graphicx}
\usepackage{indentfirst}
\usepackage{amsmath,amsfonts}
\usepackage{amstext,amsthm,amssymb}
\usepackage{mathabx}
\usepackage{mathrsfs}
\usepackage{ragged2e}
\usepackage[mathscr]{euscript}
\usepackage{multirow}
\usepackage{subfigure}
\usepackage{enumerate}
\usepackage{mathrsfs}
\usepackage{color}
\usepackage{verbatim}
\usepackage{float}

\def\und#1{\underline{#1}}
\usepackage{sectsty}

\sectionfont{\fontsize{10}{14}\selectfont}
\subsectionfont{\fontsize{10}{13}\selectfont}

\usepackage[breaklinks=true,a4paper=true,pagebackref=true,pdftex]{hyperref}
\hypersetup{
pdfauthor={Diego S. Rodrigues},
pdftitle={A semi-implicit finite element method for viscous lipid membranes}
}

\usepackage{fancyhdr}

\usepackage{authblk}

\title{\Large \bf\sc A semi-implicit finite element method for viscous lipid membranes}
\author{\small \sc Diego S. Rodrigues}
\author{\small \sc Roberto F. Ausas}
\author{\small \sc Fernando Mut}
\author{\small \sc Gustavo C. Buscaglia\thanks{gustavo.buscaglia@icmc.usp.br}}
\affil{\vspace{-0.3cm} \footnotesize \emph{Instituto de Ci\^encias Matem\'aticas e de Computa\c c\~ao,
Universidade de S\~ao Paulo, Brazil}}
\date{\relax}

\begin{document}
\maketitle
\vspace{-1.0cm}
\begin{abstract}
A finite element formulation to approximate
the behavior of lipid membranes is proposed.
The mathematical model incorporates tangential
viscous stresses and bending elastic forces, together
with the inextensibility constraint and the
enclosed volume constraint. The membrane is discretized
by a surface mesh made up of planar triangles, over
which a mixed formulation (velocity-curvature) is
built based on the viscous bilinear form (Boussinesq-Scriven operator)
and the Laplace-Beltrami identity relating position and curvature.
A semi-implicit approach is then used to discretize in
time, with piecewise linear interpolants for all variables.
Two stabilization terms are needed: The first one stabilizes
the inextensibility constraint by a pressure-gradient-projection
scheme (R. Codina and J. Blasco, Computer Methods in Applied
Mechanics and Engineering 143:373-391, 1997), the second couples curvature and velocity to
improve temporal stability, as proposed by B\"ansch (Numerische
Mathematik 88:203-235, 2001).
The volume constraint is handled by a Lagrange multiplier
(which turns out to be the internal pressure), and an
analogous strategy is used to filter out rigid-body motions.
The nodal positions are updated in a Lagrangian manner according
to the velocity solution at each time step.
An automatic remeshing strategy maintains suitable refinement
and mesh quality throughout the simulation.

Numerical experiments show the convergent and robust
behavior of the proposed method. Stability limits are
obtained from numerous relaxation tests, and convergence with
mesh refinement is confirmed both in the relaxation transient
and in the final equilibrium shape. Virtual tweezing experiments
are also reported, computing the dependence of the deformed
membrane shape with the tweezing velocity (a purely dynamical 
effect). For sufficiently high velocities, a tether develops
which shows good agreement, both in its final radius and in its
transient behavior, with available analytical solutions.
Finally, simulation results of a membrane subject to the simultaneous action of 
six tweezers illustrate the robustness of the method.
\vspace{0.2cm}

{\justifying{\bf Keywords:} Biological Membranes, Lipid Bilayer,
Canham-Helfrich Energy,
Boussinesq-Scriven Operator, Tangential Calculus, Finite Element Method.}
\end{abstract}

\section{INTRODUCTION}

Phospholipid membranes are two-molecule-thick curved surface
arrays of phospholipids \citep{alberts} that constitute the fundamental building
material of the Living Cell membrane, of many intra-cellular 
units, and of synthetic vesicles such as liposomes.
The static properties of this two-dimensional material are governed
by geometry. In fact remarkable agreement with biophysical observations
has been obtained with models in which the energy density (per unit area) is a
function of the local curvature alone \cite{canham,seifert1}. Such
an energy density is typical of \emph{elastic solids} in bending.

Numerical methods for computing equilibrium shapes of these
membranes by gradient flow (which in this context is called {\emph{Willmore
flow}}) first appeared about ten years ago, with the works
of Dziuk \cite{dziuk}, Rusu \cite{rusu05}, Feng \& Klug \cite{feng06}
and Barret {\em et al.} \cite{bgn08_jsc}, among others. These methods
evolve the geometry by gradient descent towards an equilibrium of
the applied forces (if any) with the elastic forces.
Bonito {\em et al.} \cite{bonito2010,bonito2011bulk} considered the effect of the
bulk fluid, while Elliot \& Stinner modeled two-phase effects
\cite{es10_jcp}, always in gradient flow. 

The actual dynamics of phospholipid membranes does not however
obey gradient flow. Their evolution results from the interplay 
between the applied forces, 
the hydrodynamic forces coming
from the adjacent inner and outer liquids,  and the forces that develop on
the membrane itself, which include an elastic contribution (as in
gradient flow) and also a surface viscous contribution arising
from the lipid-to-lipid sliding. 
In this article we focus just on the membrane forces, restricting 
the effect of the adjacent
liquids to just a volume constraint. The combination of the methods
proposed below with more realistic treatments of the inner and outer
liquids is straightforward (though the added computational cost is
obviously significant).  

We assume that the surface viscous forces that develop on the membrane and
determine its dynamics correspond to 
an area-preserving {\em  Newtonian surface fluid}
\citep{harland1,harland2,arroyo2}.
Our goal is thus to present a finite element method for the {\em viscous
flow} of phospholipid bilayers; i.e., for the
dynamical simulation of phospholipid bilayers, considering
an elastic model for bending deformations and a (viscous)
Newtonian area-preserving fluid model for the dissipative tangential motions. 

For this purpose, we adopt the same treatment of elastic
forces used for gradient flows \cite{dziuk,bonito2010}, combined
with a novel treatment of surface viscous forces.
The mathematical formulation of surface viscous behavior was
first derived by Scriven \cite{scriven}. Schemes for its numerical approximation
have been proposed by Arroyo and coworkers \cite{arroyo2,adh10,ra12pre}
in the axisymmetric case, and by Tasso \& Buscaglia \cite{tassobus} in the
general 3D case. The formulation of this latter article relies heavily
on the numerical differentiation of the energy of the membrane (including
an ``evanescent elasticity'' term which accounts for tangential viscosity)
to compute forces and stresses, and on yet another numerical
differentiation to compute the approximate tangent matrix. 
In this work another approach is followed, developing a semi-implicit 
scheme issued in a classical way from the continuous variational
formulation, without adjustable numerical differentiation parameters and
involving the solution of just one linear system per time step.

After introducing the mathematical formulation in Section 2 and
the proposed discretization scheme in Section 3, we assess the
proposed method through numerical examples in Section 4.
Special attention is given to experiments that involve
membrane tweezing and tether formation. 
The latter is a salient phenomenon that takes place in 
phospholipid bilayers, by which
if a small part of a vesicle is pulled away by some
localized force (using an optical trap, for example
\cite{lee})
it carries with it a narrow bilayer tube
(tether) that can be much longer that the vesicle itself and
nanometric in diameter \citep{smith_seifert}. The proposed
method is shown to be sufficiently robust to allow for accurate
simulations of tether formation and
extension, which are important to shed light on fundamental
mechanisms of cell mechanics \cite{waughI,bozic,waughII}.
Section 5 is then devoted to summarize the conclusions of
the study.

\section{MATHEMATICAL FORMULATION}

\subsection{Virtual power at the interface}

We consider the motion of a closed surface $\Gamma \subset \mathbb{R}^3$
under the action of surface elastic forces and external forces
coming from the adjacent liquid. The virtual work principle for
such a system reads
\begin{equation}
\int_\Gamma \boldsymbol{\sigma}:D_{\Gamma}{\bf v} = -d\mathcal{E}(\Gamma,{\bf v})
+\int_\Gamma {\bf f}\cdot{\bf v} \qquad \qquad \forall\,{\bf v}\,\in\,V(\Gamma),
\label{eq:virtualwork}
\end{equation}
where $\boldsymbol{\sigma}$ is the tensor of tangential stresses,
$\mathcal{E}(\Gamma)$ is the elastic energy from which elastic
forces are derived, ${\bf f}$ is the net interaction force with
the surroundings, $V(\Gamma)$ the space of admissible virtual
velocities and $D_{\Gamma}{\bf v}$ the surface virtual strain rate. 

In (\ref{eq:virtualwork}), by $d\mathcal{E}(\Gamma,{\bf v})$ we denote
the derivative (or first variation) of $\mathcal{E}(\Gamma)$ along
the virtual velocity field ${\bf v}$. In turn,  $D_{\Gamma}{\bf v}$
represents the surface differential operator
\begin{eqnarray}
  D_{\Gamma}{\bf v} = \frac12\,\mathbb{P}\,( \nabla_{\Gamma}{\bf v} + 
  \nabla_{\Gamma}{\bf v}^T )\,\mathbb{P},
\end{eqnarray}
which is the surface analog of the usual three-dimensional
symmetric gradient $D{\bf v}= (\nabla {\bf v} + \nabla {\bf v}^T)/2$.

Some elements of differential geometry are needed at this point.
We follow the presentation of Buscaglia \& Ausas \cite{ba11}, the
reader is also referred to Biria {\em et al.} \cite{bmf13} for a 
more comprehensive review.

The tensor $\mathbb{P}$ above is the tangent projector onto $\Gamma$ given by
\begin{equation}
  \mathbb{P} = \mathbb{I} - \widecheck{{\bf n}} \otimes \widecheck{{\bf n}}, 
\end{equation}
$\widecheck{{\bf n}}$ being the normal to $\Gamma$,
and the symbol $\nabla_\Gamma$ refers to the surface gradient, given by
\begin{equation}
\nabla_\Gamma f = \mathbb{P}\,\nabla \widehat{f},
\end{equation}
where $\widehat{f}$ is any smooth extension of the function $f$ 
from its values on $\Gamma$ to a three-dimensional neighborhood
of it. The surface Laplacian $\Delta_\Gamma f$ is defined as 
$\nabla_\Gamma\cdot(\nabla_\Gamma f)$.

The surface gradient $\nabla_\Gamma {\bf w}$ of a vector field ${\bf w}$
defined on $\Gamma$ is defined as the matrix (Cartesian tensor)
\begin{equation}
\{ \nabla_\Gamma {\bf w} \}_{ij}  = \{ \nabla_\Gamma w_i \}_j,
\end{equation}
where $w_i$ is the $i$-th Cartesian component of ${\bf w}$.

\subsection{The Boussinesq-Scriven operator}

The rheology of a viscous interface $\Gamma$ is 
governed by the Boussinesq-Scriven law \citep{scriven,reusken},
which is the tangential analog to the Newtonian constitutive
law, i.e.,
\begin{eqnarray}
  \boldsymbol{\sigma} = (-\,p_{\text{s}} + \lambda\,\nabla_{\Gamma} \cdot{\bf u})\,\mathbb{P}\, 
  + 2\,\mu\,D_{\Gamma}{\bf u}, \label{boussinesq}
\end{eqnarray}
where $\lambda$ and $\mu$ are surface viscosity coefficients, ${\bf u}$ is
the material velocity of the membrane particles, and $p_{\text{s}}$ is a {\em surface
thermodynamic pressure}, which requires a closure law.

An area-preserving membrane (frequently called {\em inextensible} membrane)
is defined by the constraint 
\begin{equation}
\nabla_\Gamma \cdot {\bf u}=0.
\label{eq:inextensible}
\end{equation}
The inextensible limit is obtained making $\lambda$ tend to infinity.
It is a classical result that there exists a {\em surface pressure}
$\pi_{\text{s}}$, the
Lagrange multiplier associated to the constraint (\ref{eq:inextensible}),
such that, irrespective of the closure law for $p_{\text{s}}$,
\begin{equation}
\lim_{\lambda\to +\infty} \left ( -\,p_{\text{s}} 
+\lambda \nabla_\Gamma \cdot {\bf u} \right )~=~-\pi_{\text{s}}.
\label{eq:limlambda}
\end{equation}
As a consequence, the tangential stresses from (\ref{boussinesq}) read,
for inextensible membranes,
$$
\boldsymbol{\sigma} = -\,\pi_{\text{s}} \,\mathbb{P}\, 
 + 2\,\mu\,D_{\Gamma}{\bf u}.
$$

The bilinear form that expresses the virtual power along
a virtual velocity field ${\bf v}$ performed by the stresses
$\boldsymbol{\sigma}$
corresponding to an actual velocity field ${\bf u}$ and surface
pressure $\pi_{\text{s}}$ is given by
\begin{eqnarray}
\mathscr{W}\left ( ({\bf u},\pi_{\text{s}}),{\bf v} \right ) & = &
\int_\Gamma \boldsymbol{\sigma} : D_\Gamma {\bf v} \nonumber \;\;= \\
& = & \int_\Gamma
2\,\mu\,D_\Gamma {\bf u}:D_\Gamma {\bf v} ~-~
\int_\Gamma \pi_{\text{s}}\,\nabla_\Gamma \cdot {\bf v}.
\label{eq:w}
\end{eqnarray}

{\noindent \bf Remark:} {\em
The bilinear form $\mathscr{W}$ is the surface analog of the Stokes form for
bulk fluids, namely 
$$
\mathscr{W}^{\mbox{\small bulk}}\left ( ({\bf u},p),{\bf v} \right )=
\int\,2\,\mu\,D{\bf u}:D{\bf v} - \int\,p\,\nabla\cdot {\bf v},
$$
with the integrals performed this time over the {\em volume} occupied
by the bulk fluid. As it is well known, there is a differential
operator that corresponds to $\mathscr{W}^{\mbox{\small bulk}}$, which reads
$$
-\mu\,\nabla^2 {\bf u} + \nabla p.
$$
Similarly, there exists a surface differential operator associated to
$\mathscr{W}$, which can be denoted by
$$
-\mathcal{S}_\Gamma\,{\bf u} + \nabla_\Gamma \pi_{\text{s}},
$$
but the actual expression of $\mathcal{S}_\Gamma$ is quite involved.
It can be found in the pioneering work of Scriven \cite{scriven},
which is why $\mathcal{S}_\Gamma$ is sometimes referred to as
{\em Boussinesq-Scriven operator}. It can also be found, written in
the language of differential forms, in the interesting article
by Arroyo \& DeSimone \cite{arroyo2} (see also \cite{ramos13}).
}


\subsection{The Canham-Helfrich energy}

The elastic bending energy considered here is the simplest
version of the model proposed
by Canham and Helfrich \cite{canham,helfrich},
\begin{eqnarray}
  \mathscr{E}(\Gamma) = \frac{c_{\text{\tiny{CH}}}}{2} \,\int_\Gamma {\kappa^2}, \label{CHen}
\end{eqnarray}
where $\kappa = \kappa_1 + \kappa_2$ stands for the mean scalar curvature of $\Gamma$
($\kappa_1$ and $\kappa_2$ are the principal curvatures) and
$c_{\text{\tiny{CH}}}$ 
is a material dependent parameter.
In differential geometry, equation (\ref{CHen}) is
known as Willmore energy \citep{willmore}.

The Canham-Helfrich energy (\ref{CHen}) depends on the shape of $\Gamma$
and is thus affected by motions along a virtual velocity field ${\bf v}$. 
The derivative of $\mathscr{E}$ 
along ${\bf v}$ was computed by Rusu \cite{rusu05} as
\begin{align}
d\mathscr{E}({\bf v}) \;=\;
c_{\text{\tiny{CH}}}\,\int_\Gamma \left [
\frac{|\Delta_\Gamma \boldsymbol{\chi}|^2}{2} \nabla_\Gamma \boldsymbol{\chi}
: \nabla_\Gamma {\bf v} +
\nabla_\Gamma(\Delta_\Gamma\boldsymbol{\chi}):\nabla_\Gamma {\bf v} -
2\, (\nabla_\Gamma(\Delta_\Gamma\boldsymbol{\chi})^T \widecheck{{\bf n}})
\cdot
(\nabla_\Gamma{\bf v}^T \widecheck{{\bf n}}) \right ],
\nonumber\\
\end{align}
where $\boldsymbol{\chi}$ stands for the identity mapping on $\Gamma$
(i.e., $\boldsymbol{\chi}({\bf x})={\bf x},~\forall \,{\bf x}\,\in\,\Gamma$),
which obeys
\begin{equation}
\mathbb{P}=\nabla_\Gamma \boldsymbol{\chi}, \qquad
\mbox{and}
\qquad
\boldsymbol{\kappa}\stackrel{\mbox{\tiny def}}{=}\kappa\,\widecheck{{\bf n}}
= -\Delta_\Gamma \boldsymbol{\chi}.
\end{equation}

In terms of the vector curvature $\boldsymbol{\kappa}$, the first variation
$d\mathscr{E}({\bf v})$ can be rewritten as
\begin{equation}
d\mathscr{E}({\bf v}) \;=\;
c_{\text{\tiny{CH}}}\,\int_\Gamma \left [
\frac{|\boldsymbol{\kappa}|^2}{2} \mathbb{P}
: \nabla_\Gamma {\bf v} +
(\mathbb{I}-2\,\mathbb{P}) \nabla_\Gamma {\bf v} : \nabla_\Gamma \boldsymbol{\kappa}
\right ].
\end{equation}

Equivalent formulas were produced by Dziuk \cite{dziuk} and Bonito {\em et al.} \cite{bonito2010}.
The latter was adopted in our implementation, which reads
\begin{equation}
d\mathscr{E}({\bf v}) \;=\;
c_{\text{\tiny{CH}}}\,\int_\Gamma \left [
(\mathbb{I}-2\,\mathbb{P}) \nabla_\Gamma {\bf v} : \nabla_\Gamma \boldsymbol{\kappa}
+ \frac12 (\nabla_\Gamma \cdot {\bf v})\,(\nabla_\Gamma \cdot \boldsymbol{\kappa})
\right ],
\end{equation}
which holds if $\boldsymbol{\kappa}$ obeys the weak version of 
$-\Delta_\Gamma \boldsymbol{\chi}=\boldsymbol{\kappa}$, namely
\begin{equation}
\int_\Gamma \boldsymbol{\kappa}\cdot \boldsymbol{\zeta}
=
\int_\Gamma \mathbb{P}:\nabla_\Gamma \boldsymbol{\zeta}
\qquad \qquad
\forall\,\boldsymbol{\zeta}\,\in\,H^1(\Gamma)^3.
\end{equation}

\subsection{Volume and area constraints}

Let $\mathcal{V}$ be the volume enclosed by the lipid membrane
$\Gamma$. It satisfies
\begin{equation}
\mathcal{V}=\frac13\,\int_\Gamma \boldsymbol{\chi}
\cdot \widecheck{\bf n}
\end{equation}
and its time derivative, when the membrane velocity is ${\bf u}$,
given by
\begin{equation}
\frac{d\mathcal{V}}{dt}=\int_\Gamma {\bf u}\cdot \widecheck{\bf n}.
\end{equation}

In general, osmotic equilibrium determines the (fixed) 
volume $\mathcal{V}^*$ that the surface $\Gamma$ must enclose
at all times along its evolution, so that the instantaneous
constraint reads $\int_\Gamma {\bf u}\cdot\widecheck{\bf n}=0$.
When the membrane evolution is discretized in time, however,
the enclosed volume may drift away from the value $\mathcal{V}^*$.
To mitigate this error, we implemented a volume controller
as follows
\begin{equation}
\int_\Gamma {\bf u}\cdot \widecheck{\bf n}
~=~ \frac{\mathcal{V}^*-\mathcal{V}}{\tau_v}.
\label{eq:volume}
\end{equation}
The controller drives the volume towards the target value 
$\mathcal{V}^*$ with characteristic time $\tau_v$.

Equation (\ref{eq:volume}) acts as an additional constraint on the
membrane's dynamics, which materializes as an {\em internal pressure}
$p$
(uniform) which exerts a surface force 
$$
{\bf f}_p = p\,\widecheck{\bf n}
$$
on $\Gamma$.

The area $\mathcal{A}$ of an inextensible membrane is also
constant, this time as a consequence of (\ref{eq:inextensible})
because
\begin{equation}
\frac{d\mathcal{A}}{dt}=\int_\Gamma \nabla_\Gamma \cdot {\bf u}=0.
\end{equation}
Upon time discretization, as discussed above for the enclosed volume, 
the restriction $\frac{d\mathcal{A}}{dt}=0$
may be inexactly satisfied and thus $\mathcal{A}$ may drift away from its
correct value $\mathcal{A}^*$. An area controller is thus
implemented as
\begin{equation}
\nabla_\Gamma \cdot {\bf u}  -
\frac{\mathcal{A}^*-\mathcal{A}}{\mathcal{A}\,\tau_a} = 0,
\label{eq:areacontroller}
\end{equation}
so that, integrating over $\Gamma$, one retrieves
$$
\frac{d\mathcal{A}}{dt}~=~
\frac{\mathcal{A}^*-\mathcal{A}}{\tau_a},
$$
which drives the membrane area towards $\mathcal{A}^*$ with characteristic
time $\tau_a$.

{\noindent \bf Remark:} {\em The modifications introduced by the volume and
area controllers have no effect {\em in the exact problem} if the initial 
volume equals $\mathcal{V}^*$ and the initial 
area equals $\mathcal{A}^*$. In fact, if $\mathcal{V}(t=0)=\mathcal{V}^*$
then (\ref{eq:volume}) forces $\mathcal{V}(t)$ to equal $\mathcal{V}^*$
at all times. Similarly, if $\mathcal{A}(t=0)=\mathcal{A}^*$, then
(\ref{eq:areacontroller}) implies $\mathcal{A}(t)=\mathcal{A}^*$ for all $t>0$.}

The Lagrange multiplier associated to the conservation of area is
the surface pressure $\pi_s$, already discussed, so that the area
controller adds nothing to the bilinear form (\ref{eq:w}).


\subsection{Variational formulation}

Collecting the ingredients discussed in the previous sections,
the variational formulation that determines the velocity 
of the membrane corresponds to the following {\em linear} problem:

{\noindent {\bf Problem P:} \em ``Find $({\bf u},\pi_{\text{s}},\boldsymbol{\kappa},p)\,\in\,{\bf V}\times Q \times {\bf K}
\times \mathbb{R}$
such that
\begin{eqnarray}
\hspace{1cm}\int_\Gamma
2\,\mu\,D_\Gamma {\bf u}:D_\Gamma {\bf v} \; \;-\;
\int_\Gamma \pi_{\text{s}}\,\nabla_\Gamma \cdot {\bf v} \;\;+ & & \label{eqvaru} \nonumber \\
\;\;\;+\;\;c_{\text{\tiny{CH}}}\,\int_\Gamma \left [
(\mathbb{I}-2\,\mathbb{P}) \nabla_\Gamma \boldsymbol{\kappa} : \nabla_\Gamma {\bf v}
+ \frac12 (\nabla_\Gamma \cdot \boldsymbol{\kappa})\,(\nabla_\Gamma \cdot {\bf v})
\right ] - p\,\int_\Gamma {\bf v}\cdot \widecheck{\bf n}&=& \int_\Gamma
{\bf f}\cdot {\bf v} \\
\int_\Gamma \xi\,\nabla_\Gamma \cdot {\bf u} & = & 
\frac{\mathcal{A}^*-\mathcal{A}}{\mathcal{A}\,\tau_a}\,
\,\int_\Gamma \xi \nonumber \\ & & \label{eqarea}\\
\int_\Gamma \boldsymbol{\kappa}\cdot \boldsymbol{\zeta} &=& 
\int_\Gamma \nabla_\Gamma \mathbb{P} : \nabla_\Gamma \boldsymbol{\zeta}
\nonumber \\ & & \label{eqvarkappa}\\
\int_\Gamma {\bf u}\cdot \widecheck{\bf n} & = & \frac{\mathcal{V}^*-\mathcal{V}}{\tau_v} \label{eqvolume}
\end{eqnarray}
for all $({\bf v},\xi,\boldsymbol{\zeta})\,\in\,{\bf V}\times Q \times {\bf K}$.''
}

\medskip
The surface pressure $\pi_s$, the vector curvature $\boldsymbol{\kappa}$
and the internal pressure $p$ arise in this formulation as ``by-products''
of computing ${\bf u}$. Notice that the force field ${\bf f}$ on the
right-hand side of (\ref{eqvaru}) now comprises all interaction forces
with the surroundings {\em other than that coming from the internal pressure}.
For problem {\bf P} to be well-posed,
the spaces ${\bf V}$, $Q$ and ${\bf K}$ need to be discussed. 

Assuming the surface $\Gamma$
to be smooth, which implies that $\boldsymbol{\chi}$ is smooth,
one can integrate by parts the right-hand side of (\ref{eqvarkappa})
so as to take ${\bf K}=L^2(\Gamma)^3$. There is then a unique solution
$\boldsymbol{\kappa}\,\in\,{\bf K}$, which can then be seen to 
be smooth because of the smoothness of $\boldsymbol{\chi}$. 

Let us consider then existence and uniqueness of ${\bf u}$.
For simplicity, let us set $\pi_s=p=0$ and leave aside 
Eqs. (\ref{eqarea}) and (\ref{eqvolume}),
which are constraints handled by Lagrange multipliers.
All that remains is to plug $\boldsymbol{\kappa}$
into (\ref{eqvaru}) and solve the Boussinesq-Scriven operator to
determine ${\bf u}$. 

The well-posedness of problem {\bf P} thus demands that the
bilinear form
\begin{equation}
\mathcal{B}({\bf u},{\bf v})=\int_\Gamma 2\,\mu\,D_\Gamma{\bf u}:D_\Gamma{\bf v}
\label{eqbilinear}
\end{equation}
be continuous and (weakly) coercive over the velocity space
${\bf V}$. For continuity, ${\bf V}$ must
be contained in $H^1(\Gamma)^3$. For coercivity, it must be quotiented with
the space of (infinitesimal) rigid movements
\begin{equation}
\mathcal{R}\stackrel{\mbox{\tiny def}}{=} \{
{\bf w}:\mathbb{R}^3 \to \mathbb{R}^3 ~|~
{\bf w}({\bf x})=\boldsymbol{\omega}\,\wedge\,{\bf x}\,+\,\boldsymbol{\beta},
\,\boldsymbol{\omega}, \boldsymbol{\beta}\,\in\,\mathbb{R}^3\}
\end{equation}
because $D_\Gamma {\bf w}({\bf x})=0$, for all ${\bf x}$, whenever
${\bf w}\,\in\,\mathcal{R}$. 

In this exposition we take ${\bf V}$ as {\em equal} 
to $H^1(\Gamma)^3 / \mathcal{R}$
and reason as if the bilinear form $\mathcal{B}(\bullet,\bullet )$
were coercive in ${\bf V}$. This assumption allows us to consider
${\bf u}$ as uniquely defined by (\ref{eqvaru}), assuming
$\boldsymbol{\kappa}$ already computed (and, as said, ignoring the
geometrical constraints). Problem {\bf P} is thus assumed to be well-posed,
yielding a unique solution $({\bf u},\boldsymbol{\kappa})\,\in\, {\bf V}\,\times\,{\bf K}$.

{\noindent \bf Remark:} {\em For later use, let us recall that
the energy dissipation rate of the
surface is given by
$$
\mathcal{D}=\int_{\Gamma}2\mu\,\|D_\Gamma{\bf u}\|^2=\mathcal{B}({\bf u},{\bf u}).
$$
}

If we now consider the inextensibility equation (\ref{eqarea}), the
situation is similar to that of the incompressible Stokes equation
in that an inf-sup condition arises, namely,
\begin{equation}
\inf_{0\neq\xi\,\in\,Q}\sup_{0\neq\bf v\,\in\,{\bf V}} 
\frac{\int_\Gamma \xi\,\nabla_\Gamma \cdot {\bf v}}
{\|\xi\|_Q\,\|{\bf v}\|_{\bf V}} ~>~0.
\end{equation}
We assume that this condition is fulfilled when $Q=L^2(\Gamma)$.

\medskip

The reader should be warned that the viscous model above does not
incorporate the layer-to-layer slippage of the two molecular sheets
that form the lipid bilayer. This mode of deformation may well be
dominant in some situations, as discussed by Evans \& Yeung \cite{ey94cpl}
and more recently by Rahimi \& Arroyo \cite{ra12pre}. In this 
contribution the focus is in the numerical treatment of
the Boussinesq-Scriven operator coupled to the Canham-Helfrich
elastic model, so that the incorporation of layer-to-layer slippage
models is left for future work.

\subsection{The evolutionary problem}

Up to now we have considered a single instant of time, at which the membrane
configuration is described by a surface $\Gamma$. Since an outcome of
the instantaneous problem is in fact the velocity field with which the
membrane's particles are moving, one is lead to the
following evolutionary problem:

{\noindent \bf Evolutionary problem EP:} {\em ``Given $\Gamma(0)$, the initial surface,
compute the continuous family of surfaces $\Gamma(t)$ that evolves from $\Gamma(0)$
as convected by
the velocity field $\bf{u}(t):\Gamma(t)\to \mathbb{R}^3$
that solves problem P. In mathematical terms, the
family $\Gamma(t)$ must satisfy
\begin{equation}
\forall {\bf x}\in\Gamma(t), ~
\mbox{dist}\left ({\bf x}+{\bf u}({\bf x},t)\,\delta\,t, \Gamma(t+\delta t) \right )
~\leq C\,\delta t^2 ,
\label{kinematic}
\end{equation}
where $\mbox{dist}$ stands for the distance between a point
and a surface, for some $C>0$.
}

Notice that the tangential component of ${\bf u}(t)$ is inconsequential
in the evolution of $\Gamma(t)$. However, and contrary to what happens
in gradient flow, the tangential velocity generated by viscous flow 
is {\em not} zero.

\section{DISCRETIZATION}

We consider triangulation surfaces in 3D space, which for a fixed
mesh connectivity are uniquely described by the vector $\und{\bf X}$ of vertex
positions. Time is discretized so that a sequence of triangulation
surfaces $\Gamma^0$, $\Gamma^1$,$\ldots$,$\Gamma^n$,$\ldots$ 
are computed, corresponding to vertex positions 
$\und{\bf X}^0$, $\und{\bf X}^1$,$\ldots$,$\und{\bf X}^n$,$\ldots$.

On each $\Gamma^n$ we define the piecewise-affine finite element space
\begin{equation}
\mathbb{P}^n_1 = \{f\,\in\,\mathcal{C}^0(\Gamma^n)
~:~f|_K~\mbox{is affine},
~\forall\,K~\mbox{triangle in}~\Gamma^n \}
\end{equation}
and the approximation spaces for velocity, surface pressure and
curvature
\begin{eqnarray}
{\bf V}_h^n&=&\left ( \mathbb{P}^n_1 \right )^3\,/\,\mathcal{R}, \\
{Q}_h^n&=& \mathbb{P}^n_1, \\ 
{\bf K}_h^n&=&\left ( \mathbb{P}^n_1 \right )^3. 
\end{eqnarray}

{\noindent \bf DISCRETE PROBLEM DP:}
Defining $\delta t = t_{n+1}-t_n$, the proposed scheme 
updates the nodal positions in a Lagrangian way, i.e.,
\begin{equation}
{\bf X}^{J,n+1}={\bf X}^{J,n}+\delta t\,{\bf u}_h^{n+1}({\bf X}^{J,n}),
\label{equpdate}
\end{equation}
where $J$ is the nodal index,
so that (\ref{kinematic}) is by construction satisfied.
Notice that the velocity field ${\bf u}_h^{n+1}$ is computed
on $\Gamma^n$ and is thus an element of ${\bf V}_h^n$. 

The fully discrete linear problem that determines ${\bf u}_h^{n+1}$ is
the following:

\medskip

{\noindent \em ``Find $({\bf u}_h^{n+1},\pi_h^{n+1},\boldsymbol{\kappa}_h^{n+1},p^{n+1})\,\in\,{\bf V}_h^n\times Q_h^n \times {\bf K}_h^n \times \mathbb{R}$
such that
\begin{eqnarray}
\hspace{1cm}\int_{\Gamma^n}
2\,\mu\,D_\Gamma {\bf u}_h^{n+1}:D_\Gamma {\bf v} \;+
\int_{\Gamma^n} \pi_h^{n+1}\,\nabla_\Gamma \cdot {\bf v} \;\;-\,p^{n+1}\,\int_{\Gamma^n} {\bf v}\cdot \widecheck{\bf n}\,+ & & \label{eqvaruh} \nonumber \\
\;\;\;+\;\;c_{\text{\tiny{CH}}}\,\int_{\Gamma^n} \left [
(\mathbb{I}-2\,\mathbb{P}) \nabla_\Gamma {\bf v} : \nabla_\Gamma \boldsymbol{\kappa}_h^{n+1}
+ \frac12 (\nabla_\Gamma \cdot {\bf v})\,(\nabla_\Gamma \cdot \boldsymbol{\kappa}_h^{n+1})
\right ] & = &\int_{\Gamma^n} {\bf f}^{n+1}\cdot{\bf v}
\\
\int_{\Gamma^n} \xi\,\nabla_\Gamma \cdot {\bf u}_h^{n+1}
+ \int_{\Gamma^n} \gamma_h\,(\nabla_\Gamma \pi_h^{n+1}-{\bf g}_h^n)\cdot
\nabla_\Gamma \xi
 & = & 
\frac{\mathcal{A}^*-\mathcal{A}^n}{\mathcal{A}^n\,\tau_a}\,
\,\int_{\Gamma^n} \xi
\nonumber \\
& &
\label{eqvarpih}
\\
-\,\int_{\Gamma^n} \tau_\kappa \nabla_\Gamma {\bf u}_h^{n+1}:\nabla_\Gamma 
\boldsymbol{\zeta}+
\int_{\Gamma^n} \boldsymbol{\kappa}_h^{n+1}\cdot \boldsymbol{\zeta} &=& 
\int_{\Gamma^n}  \mathbb{P}: \nabla_\Gamma \boldsymbol{\zeta}\hspace{1cm}
\label{eqvarkappah} \\
\int_{\Gamma^n} {\bf u}_h^{n+1}\cdot \widecheck{\bf n} & = & \frac{\mathcal{V}^*-\mathcal{V}^n}{\tau_v} \label{eqvolumeh}
\end{eqnarray}
hold $\forall {\bf v} \in{\bf V}_h^n$, $\forall \xi\in Q_h^n$ and
$\forall \boldsymbol{\zeta}\in {\bf K}_h^n$.''}
Together with (\ref{equpdate}), this completely defines the fully
discrete formulation. Notice that all integrals are performed over the
known discrete surface $\Gamma^n$.

Several remarks are in order:

\begin{itemize}
\item Algorithms that compute the velocity with frozen vertex positions,
as is the case of {\bf DP}, suffer severe stability restrictions on $\delta t$.
The trend has thus been to ``implicitize'' as many terms as possible
while keeping the problem to be solved at each time step linear, 
as done by Rusu \cite{rusu05}, Dziuk \cite{dziuk} and others. 
\item A stabilization term $$\int_{\Gamma^n} \gamma_h\,(\nabla_\Gamma \pi_h^{n+1}-{\bf g}_h^n)\cdot
\nabla_\Gamma \xi$$
has been added in the inextensibility equation (\ref{eqvarpih}). This
aims at stabilizing checkerboard modes arising from the equal-order
interpolation of ${\bf u}_h$ and $\pi_h$. The stabilization technique
is taken from the ``stabilization by pressure gradient projection'' method
proposed by Codina \& Blasco \cite{cb97,bbc00,cbbh01}. The vector field
${\bf g}^n$ is the $L^2(\Gamma)$-projection of $\nabla_\Gamma \pi_h^n$ onto
$(Q_h^n)^3$, i.e.,
\begin{equation}
\int_{\Gamma^n} {\bf g}_h^n\cdot {\bf v} = 
\int_{\Gamma^n} \nabla_\Gamma \pi_h^n \cdot {\bf v}
\qquad
\qquad \forall \,{\bf v}\,\in\,(Q_h^n)^3.
\label{eqgh}
\end{equation}
The parameter $\gamma_h$ varies from element to element,
according to
\begin{equation}
\gamma_h = \frac{h_K^2}{10\,\mu},
\label{eqgammah}
\end{equation}
where $h_K$ is the diameter of element $K$. The consistent mass matrix
is used in solving (\ref{eqgh}).
\item By comparing (\ref{eqvarkappah}) to its exact version
(\ref{eqvarkappa}), one notices the addition of the
stabilization term due to B\"ansch \cite{bansch01}
$$
-\,\int_{\Gamma^n} \tau_\kappa \nabla_\Gamma {\bf u}_h^{n+1}:\nabla_\Gamma 
\boldsymbol{\zeta}_h,
$$
for which the usual choice is $\tau_\kappa = \delta t$, adopted throughout
this article. This term significantly increases the temporal
stability. It allows time steps hundreds of times larger than 
those allowed by the unstabilized algorithm ($\tau_K=0$).
\item The space ${\bf V}_h^n$ needs to have its rigid modes
filtered out. We accomplish this by a classical Lagrange multiplier
technique, which adds 6 equations ($\int_{\Gamma^n}{\bf u}_h^{n+1}={\bf 0}$
and $\int_{\Gamma^n}{\bf x}\wedge{\bf u}_h^{n+1}={\bf 0}$) and 6 unknowns to the
global matrix.
\item The characteristic times $\tau_a$ and $\tau_v$ of the area and
volume controllers, respectively, which are non-physical, are taken as
\begin{equation}
\tau_a~=~\tau_v~=~10\,\delta t.
\label{eqtauav}
\end{equation}
This choice yields the best results in terms of accuracy and 
stability, as concluded from numerous experiments. 
\end{itemize}

\section{REMESHING}

The simulation of evolving surfaces that undergo large deformations requires
adaptive meshing techniques to mantain good accuracy
along the computations. The loss of acuracy is not only related to the
degradation of triangles quality, but also to the changes in time of
the local surface curvature.
In order to cope with these issues an automatic discrete surface
regridding software was employed \cite{Lohner19961}.

The remeshing procedure starts by defining a single discrete patch as
the whole support surface whose boundary is the largest edge. Using this
edge as the initial front, the discrete patch is triangulated using an
advancing front technique. The desired local element size is defined using the
curvature information by the rule 
$$h^*(\kappa) = \frac{c_h}{\kappa}, $$
where
$c_h$ is a user-defined parameter, and $\kappa$ is the scalar curvature
provided by the field solver. The specified
element size is isotropic since only scalar curvature information is used.
The output of this step is a completely new discrete surface. Although the
new nodes lie on the original surface, the two surfaces are not 
coincident. In particular, discrepancies in the
curvature introduce discontinuities in the elastic energy after 
each remeshing. 
These perturbations are however rapidly dissipated and seem to 
not have any major impact on the simulation results.

In order to assess the quality of a given surface discretization, two
parameters are defined as measures of the shape and size quality of
each individual triangle $K$ as follows:
\begin{itemize}
\item Element {\em shape} quality: 
$$q_K^{\mbox{\tiny{shape}}} = (12\,\sqrt{3})\,A_K / P_K^2,$$ 
where $A_K$
and $P_K$ are the triangle's area and perimeter, respectively.
\item Element {\em size} quality:
$$
q_K^{\mbox{\tiny{size}}} = \min \left\{\frac{h^*(\kappa_K)}{h_K},1\right\}.
$$ 
\end{itemize}
Global measures of shape and
size qualities are then defined as
$$
Q_{\mbox{\tiny{shape}}} = \min_K \left\{q_K^{\mbox{\tiny{shape}}}\right\}
\qquad \mbox{and}\qquad
Q_{\mbox{\tiny{size}}} = \min_K \left\{q_K^{\mbox{\tiny{size}}}\right\},
$$
respectively.

The evolving discrete surface is remeshed every time one
of the two quality measures drops below given threshold values
($Q_{\mbox{\tiny{shape}}}^*$ and $Q_{\mbox{\tiny{size}}}^*$).
For all the simulations presented in this paper
$Q_{\mbox{\tiny{shape}}}^*$ and $Q_{\mbox{\tiny{size}}}^*$
were set to $0.65$ and $0.55$ respectively.

\section{NUMERICAL RESULTS}

\subsection{Adimensionalization}

It is convenient to express the numerical results in non-dimensional
form. For this purpose, one defines the basic length scale
for an inextensible membrane as 
$$
R_0 = \sqrt{\frac{\mathcal{A}}{4\,\pi}},
$$
so that the non-dimensional area is always $4\,\pi$.
This allows for the definition of consistent scales for
velocity, surface pressure, surface stress, curvature, internal pressure
and other variables as shown in Table \ref{tabdimensions}.

\begin{table}[hp]
\begin{center}
\begin{tabular}{|l|c|c|c|}
\hline
& & & \\
{\bf Quantity} & {\bf Symbol} & {\bf Scale} & {\bf Sample value}\\
& & & \\
\hline
Space & ${\bf x}$ & $R_0$ & $10^{-6}$ m \\
\hline
Time &$t$& $\dfrac{\mu\,R_0^2}{c_{\text{\tiny{CH}}}}$ & $0.25$ s \\
\hline
Velocity &${\bf u}$& $\dfrac{c_{\text{\tiny{CH}}}}{\mu\,R_0}$ & $4\times 10^{-6}$ m/s \\
\hline
Area & $\mathcal{A}$ & $R_0^2$ & $10^{-12}$ m$^2$ \\
\hline
Energy & $\mathcal{E}$ & $c_{\text{\tiny{CH}}}$ & $4\times 10^{-20}$ J \\
\hline
Dissipation & $\mathcal{D}$ &$\dfrac{c^2_{\text{\tiny{CH}}}}{\mu\,R^2_0}$ & $1.6\times 10^{-19}$ W \\
\hline
Surface pressure & $\pi_s$ & $\dfrac{c_{\text{\tiny{CH}}}}{R_0^2}$ & $4\times 10^{-8}$ Pa-m \\
\hline
Surface stress & $\boldsymbol{\sigma}$ & $\dfrac{c_{\text{\tiny{CH}}}}{R_0^2}$ & $4\times 10^{-8}$ Pa-m \\
\hline
Curvature & $\boldsymbol{\kappa}$ & $\dfrac{1}{R_0}$ & $10^{6}$ m$^{-1}$ \\
\hline
Internal pressure & $p$ & $\dfrac{c_{\text{\tiny{CH}}}}{R_0^3}$ & $0.04$ Pa \\
\hline
Surface force & ${\bf f}$ & $\dfrac{c_{\text{\tiny{CH}}}}{R_0^3}$ & $0.04$ Pa \\
\hline
Force & ${F}$ & $\dfrac{c_{\text{\tiny{CH}}}}{R_0}$ & $4\times 10^{-14}$ N \\
\hline
\end{tabular}
\end{center}
\caption{Adimensionalization scales for the intervening variables.
The sample values correspond to $R_0=10^{-6}$ m, 
$c_{\text{\tiny{CH}}}=4\times 10^{-20}$ J and $\mu = 10^{-8}$ Pa-s-m.}
\label{tabdimensions}
\end{table}

A {\em relaxation experiment} corresponds to solving problem {\bf P}
repeatedly starting from an initial configuration $\Gamma^0$ and with
no forces other than the internal pressure applied (i.e.; ${\bf f}=0$),
so that the membrane evolves towards a nearby equilibrium.
In a relaxation experiment all non-dimensional variables depend just
on the non-dimensional initial configuration $\widehat{\Gamma}^0$,
where $\widehat{\Gamma}^0$ is the scaled version of $\Gamma^0$, i.e.;
\begin{equation}
{\bf x}\,\in\,\Gamma^0 ~\Leftrightarrow~
\widehat{\bf x}\stackrel{\mbox{\tiny def}}{=}
\frac{1}{R_0}\,{\bf x}~\,\in\,\widehat{\Gamma}^0.
\end{equation}
If two relaxation experiments share the same $\widehat{\Gamma}^0$,
then the time history (in terms of non-dimensional time) of
all (non-dimensional) variables must coincide, irrespective of
the actual values of $R_0$, $c_{\text{\tiny{CH}}}$ and $\mu$.

In a {\em tweezing experiment}, on the other hand, there is a part of the
membrane that is pulled away with some imposed velocity $V_T$ or 
some imposed force $F_T$. In this case the non-dimensional
solutions will depend both on $\widehat{\Gamma}^0$
and on the non-dimensional value of the imposed velocity
or force, which acts as an additional non-dimensional parameter.

In what follows, all reported quantities are non-dimensional
unless explicitly said otherwise. The sample values tabulated above
may help the reader in translating the non-dimensional results
into physical quantities.

\subsection{Relaxation experiments: Stability limit, convergence and 
equilibrium shape}

Equilibrium shapes of lipid membranes, or equivalently 
stationary points (local minima) of the Canham-Helfrich energy, are
configurations $\Gamma^\infty$ at which the membrane is in static equilibrium
(the solution to problem {\bf P} is ${\bf u}({\bf x})=0\,\forall
{\bf x}\,\in\,\Gamma$). Equilibrium shapes have been studied
extensively by Seifert and coworkers \cite{seifert1}, among others.

The {\em viscous relaxation} of a membrane corresponds to the
evolution, without any external force (${\bf f}\equiv 0$),
from an initial shape $\Gamma^0$ towards an equilibrium
shape $\Gamma^\infty$, obeying the viscous model described in this
article. In what follows we assess the performance of the proposed
method (defined by Eqs. (\ref{equpdate})-(\ref{eqvolumeh})) for
relaxation experiments. For this purpose, we first determine 
the stability limit of the method (maximum $\delta t$ for
stable behavior, as a function of the mesh
size), and then conduct numerical relaxations
with increasingly refined meshes. There is
no analytical solution for the relaxation transient, so that 
what is being analyzed is the consistency of the results obtained
for different meshes. The discrete equilibrium shape, on the
other hand, can be compared to quasi-analytical results (analogous
to those of Veerapaneni {\em et al.} \cite{veerapaneni}).

\subsubsection{Stability limits}

\begin{figure}
\begin{center}
\scalebox{1.4}{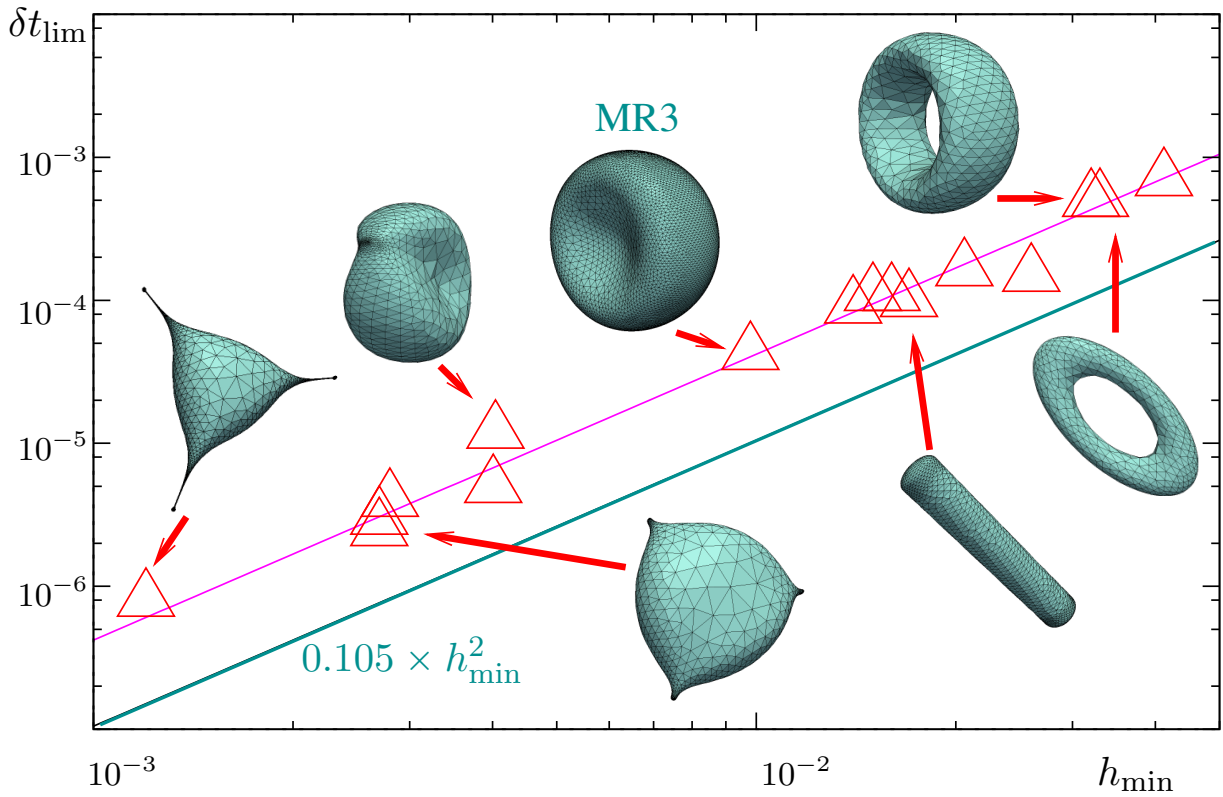}
\end{center}
\caption{Stability limit $\delta t_{\lim}$ plotted as a function of the
minimum edge size $h_{\min}$. The triangles are experimentally obtained
values for meshes of different shapes and refinement (some of the shapes
are shown and the corresponding data point indicated). In magenta 
the best-fit line $0.42\,h_{\min}^2$. In cyan the adopted time-stepping
strategy, $0.105\,h_{\min}^2$.
}
\label{figdtliminserts}
\end{figure}

The initial shape can be seen in Figure \ref{figdtliminserts},
with a triangulation that corresponds to the finest mesh employed
(mesh MR3).
The enclosed volume is $\mathcal{V}(t=0)=3.1907$, and this same
value is taken as $\mathcal{V}^*$. Though this value is
non-dimensional, it is customary to express the volume in terms
of another non-dimensional quantity, the {\em reduced volume} \cite{seifert1}
\begin{equation}
v~\stackrel{\mbox{\tiny def}}{=}~\frac{6\,\sqrt{\pi}\,\mbox{Volume}}
{\mbox{Area}^{\frac32}}~=~\frac{3\,\mathcal{V}}{4\,\pi},
\end{equation}
where ``Volume'' and ``Area'' stand for the actual (dimensional) volume enclosed
by the membrane and area of the membrane, respectively. 
The reduced volume enclosed by mesh MR3 is $v(t=0)=0.7617$.

All results below and in the next sections are computed with
algorithm DP (Equations (\ref{equpdate})-(\ref{eqvolumeh})),
with $\tau_K=\delta t$, $\tau_a=\tau_v=10\,\delta t$ and
$\gamma_h$ given by (\ref{eqgammah}).

The first experiments aim at determining the maximum time step
size $\delta t_{\lim}$ for which the fully-discrete method DP
behaves in a stable way.
For this purpose, one hundred time steps are run on each
mesh for several choices of $\delta t$. Unstable runs
are easily recognizable by violent fluctuations of the elastic
energy and of the maximum velocity. The limit value $\delta t_{\lim}$
is obtained by dychotomic search with a tolerance $\leq 20\%$.

Three increasingly refined triangulations are employed, of which
the most refined is the already described mesh MR3.
The maximum time steps allowed by the method can be observed in 
Table \ref{tabledtlim}. They obey the formula
\begin{equation}
\delta t_{\lim}~\simeq~0.42\,h_{\min}^2
\label{eqdtlim}
\end{equation}
almost exactly. Notice that this formula is non-dimensional,
expressed dimensionally it reads
$$
\delta t_{\lim}~\simeq~\frac{0.42\,\mu}{c_{\mbox{\tiny{CH}}}}\,h_{\min}^2
\qquad
\mbox{(dimensionally).}
$$
The constant 0.42 can of course depend on the shape of the membrane, so
that a similar study was performed on several very different shapes and
with uniform or adaptively refined triangulations. The $\delta t_{\lim}$
obtained for each initial mesh is plotted as a function of $h_{\min}$
in Figure \ref{figdtliminserts}. 

The best-fit line in magenta corresponds
to (\ref{eqdtlim}), which as observed from the plot in some cases
overestimates $\delta t_{\lim}$. Further, we have observed quite often
that choosing $\delta t$ very close to the stability limit deteriorates
the accuracy of the computations. This could be a consequence of the term
$\int_{\Gamma^n} \tau_K \nabla_\Gamma {\bf u}_h^{n+1}:\nabla_\Gamma \boldsymbol{\zeta}$
in (\ref{eqvarkappah}), since we are taking $\tau_K=\delta t$. For these
two reasons we adopt as automatic time-step determination formula ({\em adjusted
every single time step}) one fourth of the value given by (\ref{eqdtlim}),
that is,
\begin{equation}
\delta t~=~\delta t^*(h_{\min})~\stackrel{\mbox{\tiny{def}}}{=}~0.105\,h_{\min}^2.
\label{eqdeltat}
\end{equation} 
Unless otherwise stated, all relaxation experiments described below have
been conducted with this time-stepping strategy.

\begin{table}
\begin{center}
\begin{tabular}{|l|c|c|c|c|}
\hline
& & & & \\
{\bf Mesh} & \# nodes & \# elements & $h_{\min}$ & $\delta t_{\lim}$\\
& & & & \\
\hline
& & & & \\
{\bf MR1} & 592 & 1180 & 0.041 & $7.0\times 10^{-4}$\\
& & & & \\
\hline
& & & & \\
{\bf MR2} & 2177 & 4350 & 0.021 & $1.6\times 10^{-4}$\\
& & & & \\
\hline
& & & & \\
{\bf MR3} & 8126 & 16248 & 0.010 & $4.2\times 10^{-5}$\\
& & & & \\
\hline
\end{tabular}
\end{center}
\caption{Maximum time step for stable behavior of the method, as
obtained for each of the meshes of the relaxation study.}
\label{tabledtlim}
\end{table}

\subsubsection{Convergence of relaxation dynamics}

Let us assess now the convergence of the proposed method.
The initial meshes are MR1, MR2 and MR3, and the time step is
updated according to (\ref{eqdeltat}). The initial values
of the time step are, thus, $1.75\times 10^{-4}$, 
$0.40 \times 10^{-4}$ and $1.05\times 10^{-5}$. The simulated
non-dimensional time is $0.06$.

There is no analytical
solution for this evolutionary problem, so that the experiments
aim at checking the consistency of the results with mesh
(and time step) refinement.

\begin{figure}[htp]
\begin{center}
\scalebox{0.7}{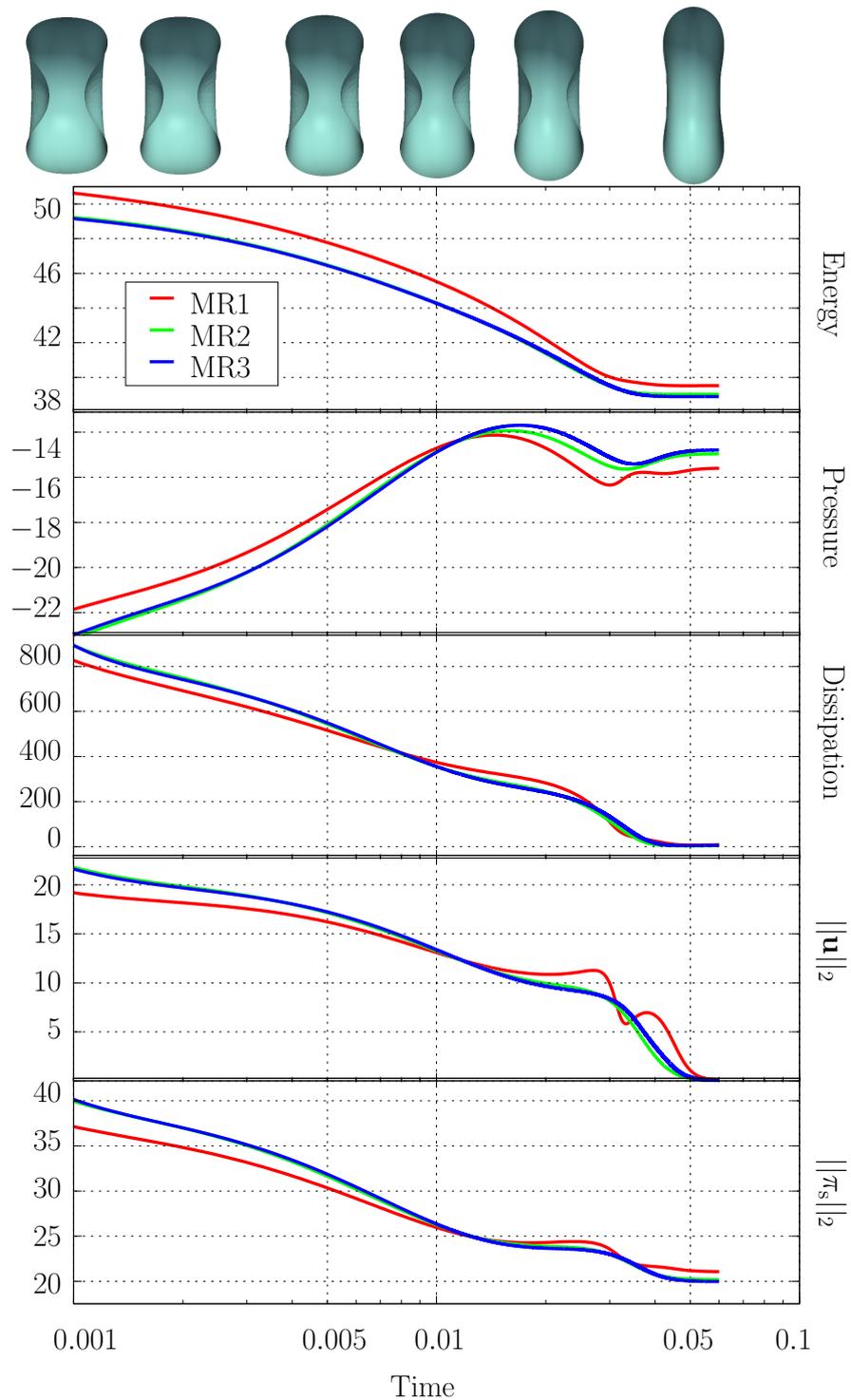}
\end{center}
\caption{
Time evolution of energy, internal pressure, dissipation, velocity
norm and surface pressure norm (both in $L^2(\Gamma)$) along the
relaxation experiment. The different colors correspond to the
increasingly refined meshes MR1, MR2 and MR3. On top, the shape
of the membrane at different instants (the horizontal position of
each shape approximately corresponds to its time).
}
\label{figrelaxation}
\end{figure}

In Figure \ref{figrelaxation} we plot several integral
quantities of the relaxation process, namely the
(all non-dimensional) energy $\mathcal{E}$, the internal
pressure $p$, the dissipation rate $\mathcal{D}$,
and the $L^2(\Gamma)$ norms of the velocity ${\bf u}$
and of the surface pressure $\pi_s$, as functions of non-dimensional
time $t$. The shape evolution is shown at the top of the Figure.

The relaxation process is seen to take until
about $t=0.06$, with an energy reduction of about 20\%
(from $\sim 49$ to $\sim 39$). The consistency of the
curves corresponding to MR1, MR2 and MR3, and the
close agreement between the two finest meshes, provide 
strong evidence of mesh convergence. Notice how some
spurious transient that takes place at $t\simeq 0.03-0.04$
for mesh MR1 (especially evident in the plot of 
$\|{\bf u}\|_2$) completely disappears after mesh refinement.

\subsubsection{Equilibrium shapes}

Discrete equilibrium shapes can be obtained by gradient flow
or by viscous flow, once the evolutionary problem reaches its
steady state. It is an important consistency check for the
proposed method that the discrete equilibrium shapes it provides
are indeed approximations of exact equilibrium shapes.

To perform this check, quasi-analytical solutions were computed
for axisymmetric shapes by numerically integrating the associated
system of ODEs with an extremely fine discretization. 
In this way, axisymmetrical versions of the
quasi-analytical shapes produced by Veerapaneni {\em et al.} 
\cite{veerapaneni} were obtained. They can be compared
to the numerical shapes at which the algorithm arrives after the
relaxation process. We selected for this comparison oblate
equilibrium shapes with reduced volumes of $v=0.61$ and
$v=0.81$. For each reduced volume, three increasingly refined
meshes were used, as in the previous section (in fact, essentially
the same meshes).

To compare the 3D results with the axisymmetric solution,
the symmetry axis of the 3D mesh is identified by diagonalizing
the tensor $\int_\Gamma {\bf x}\otimes {\bf x}$, and so a
cylindrical coordinate system $r-z-\phi$ can be assigned to 
each point in $\Gamma$, and also an arc-length coordinate $s$
along the meridians. 

\begin{figure}[htp]
\begin{center}
\scalebox{0.49}{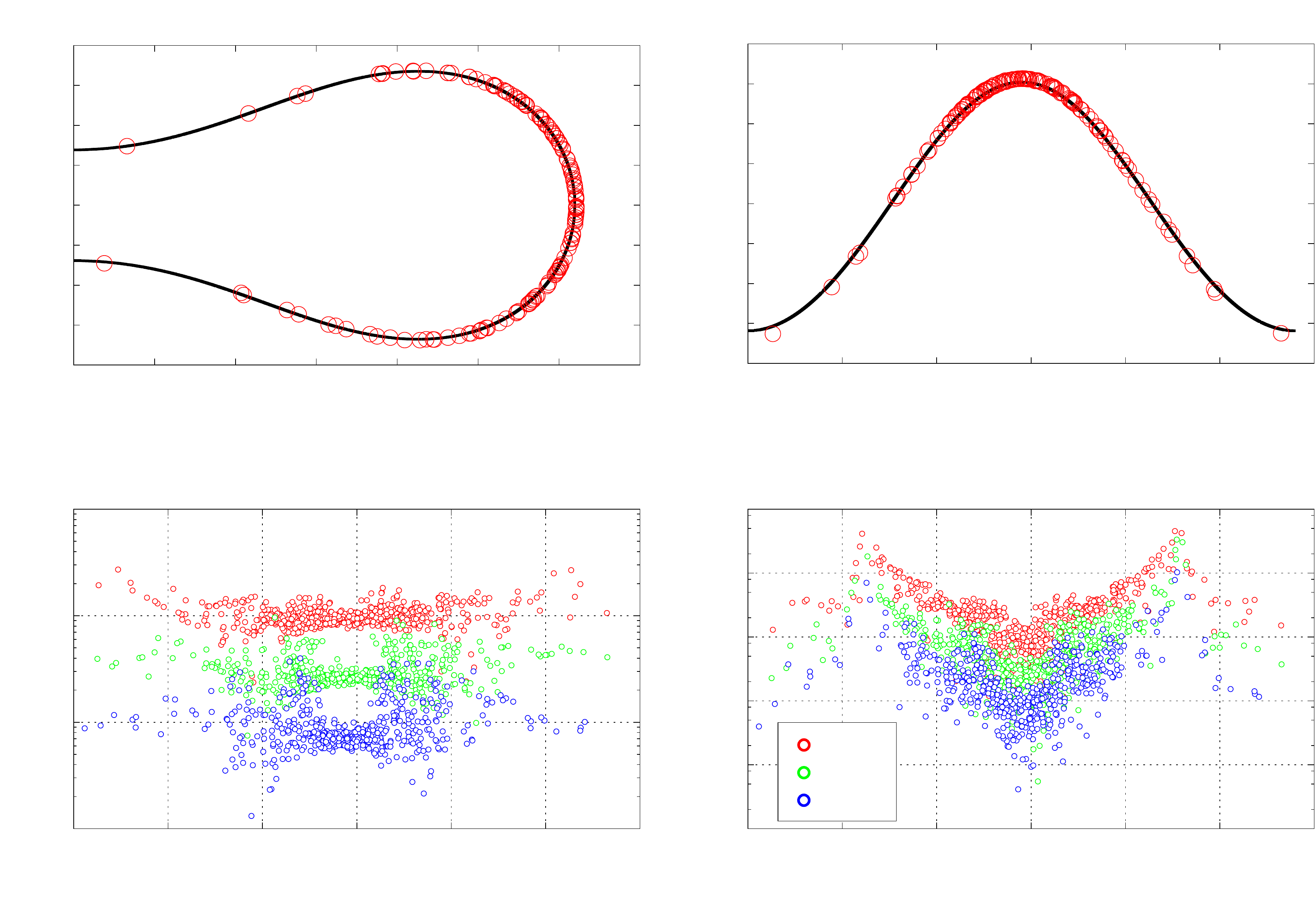}
\end{center}
\caption{Numerical and quasi-analytical equilibrium shapes for $v=0.61$.
(a) Coordinates $r-z$ of nodal positions (red circles) and exact
shape (black line). (b) Mean curvature at the nodes as a function of the
meridian-arc length (red circles) and exact mean curvature (black line).
(c) Nodal errors of the curvature vector $\boldsymbol{\kappa}_h$ as
a function of the meridian-arc length for the three meshes MR1 (red), 
MR2 (green) and MR3 (blue). (d) Idem as (c) for the normal vector
$\widecheck{\bf n}_h$.}
 \label{figequilibrium}
\end{figure}

In Figure \ref{figequilibrium}(a) the
$r-z$ coordinates of the nodes of the {\em coarsest} relaxed mesh
are superposed to the corresponding quasi-analytical curves (just
one fifth of the nodes are plotted, to leave the exact curve visible). 
Just the results of the equilibrium shape corresponding to $v=0.61$ are shown,
since those of $v=0.81$ are analogous.
The shape is seen to be quite correctly reproduced. To compare
the curvature, we plot it as a function of the arc-length
coordinate in Figure \ref{figequilibrium}(b). Each data point
of these figures involves an error, from which we compute
$$
\mbox{err}({\bf x})=\left [\frac{1}{\mbox{\# nodes}}
\sum_{J\,\in\,\mbox{\tiny{nodes}}} \|{\bf X}^J - \texttt{cp}({\bf X}^J)\|^2
\right ]^{\frac12},
$$
where \texttt{cp}$({\bf x})$ is the closest projection of ${\bf x}$
onto the exact equilibrium shape $\Gamma$. In the same way,
comparing the numerical nodal values of the different quantities
to their exact value at the closest point of $\Gamma$,
we compute discrete estimates of the errors of the different
fields, i.e., $\mbox{err}(\widecheck{\bf n})$, 
$\mbox{err}(\boldsymbol{\kappa})$, $\mbox{err}(\kappa)$, 
 $\mbox{err}(\pi_s)$. 

The results are summarized in Table \ref{tabequilibriumshape} for
the two reduced volumes $v=0.61$ and $v=0.81$. One observes convergent
behaviors of order $\mathcal{O}(h^{5/3})$ for the position
and $\mathcal{O}(h^{3/2})$ for the vector curvature $\boldsymbol{\kappa}$,
which are the main unknowns of the problem. The surface pressure
$\pi_s$ seems to converge with first order, while internal pressure $p$
and the elastic energy $\mathcal{E}$ seem to be second order.

It is interesting that the algorithmic normal $\widecheck{\bf n}$,
numerically computed as $\boldsymbol{\kappa}_h/\kappa_h$, converges
with less accuracy than $\boldsymbol{\kappa}$ itself. In Figure
\ref{figequilibrium}(c)-(d) we plot the error distribution of $\boldsymbol{\kappa}$
and of $\widecheck{\bf n}$ as a function of the arc-length coordinate.
Notice how the error in $\widecheck{\bf n}$ concentrates at regions
where the mean curvature takes values close to zero.

\begin{table}
{\small
\centerline{{\bf Case} $v=0.6104$ ($p=-15.61$, $\mathcal{E}=48.47$)}
\medskip
\begin{center}
\begin{tabular}{ccccccccc}
\hline
& & & & & & & & \\
mesh &  $h_{\min}$  &  $\mbox{err}({\bf x})$    &    
$\mbox{err}(\widecheck{\bf n})$ &
$\mbox{err}(\boldsymbol{\kappa})$ &
$\mbox{err}(\kappa)$ &
$\mbox{err}(\pi_s)$ &
$\mbox{err}(p)$ &
$\mbox{err}(\mathcal{E})$ \\ 
& & & & & & & & \\
\hline
& & & & & & & & \\
MR1  & 0.04 & 2.91e-03 & 1.21e-01 & 1.04e-01 & 7.34e-02 & 6.26e-01 & 5.1581e-01 & 5.7932e-01 \\
& & & & & & & & \\
MR2  & 0.02 & 8.45e-04 & 3.43e-02 & 3.25e-02 & 1.99e-02 & 2.34e-01 & 1.4022e-01 & 1.5458e-01 \\
& & & & & & & & \\
MR3  & 0.01 & 2.68e-04 & 2.59e-02 & 1.23e-02 & 5.38e-03 & 1.55e-01 & 3.7460e-02 & 3.8833e-02 \\
& & & & & & & & \\
EOC   &       & 1.67 & 1.08 & 1.50 & 1.83 & 0.98 & 1.84 & 1.90 \\ 
& & & & & & & & \\
\hline
\end{tabular}

\medskip
\centerline{{\bf Case} $v=0.8101$ ($p=-14.39$, $\mathcal{E}=35.89$)}
\medskip
\begin{tabular}{ccccccccc}
\hline
& & & & & & & & \\
mesh &  $h_{\min}$  &  $\mbox{err}({\bf x})$    &    
$\mbox{err}(\widecheck{\bf n})$ &
$\mbox{err}(\boldsymbol{\kappa})$ &
$\mbox{err}(\kappa)$ &
$\mbox{err}(\pi_s)$ &
$\mbox{err}(p)$ &
$\mbox{err}(\mathcal{E})$ \\ 
& & & & & & & & \\
\hline
& & & & & & & & \\
MR1  &  0.04 & 2.61e-03 & 1.48e-01 & 5.41e-02 & 3.90e-02 & 5.74e-01 & 3.92e-01 & 1.49e-01  \\
& & & & & & & & \\
MR2  & 0.02 & 7.54e-04 & 4.43e-02 & 1.62e-02 & 1.05e-02 & 1.74e-01 & 1.08e-01 & 3.52e-02 \\
& & & & & & & & \\
MR3  & 0.01 & 2.55e-04 & 2.52e-02 & 6.02e-03 & 2.80e-03 & 8.12e-02 & 2.91e-02 & 5.36e-03 \\
& & & & & & & & \\
EOC   &       & 1.63 & 1.24 & 1.54 & 1.85 & 1.37 & 1.83 & 2.34 \\ 
& & & & & & & & \\
\hline
\end{tabular}
\end{center}
\caption{Experimental convergence analysis of the different variables
as compared to those of the exact shapes for $v=0.61$ and $v=0.81$.
EOC stands for ``estimated order of convergence''.}
\label{tabequilibriumshape}
}
\end{table}

\subsubsection{Mesh distortions near equilibrium}

The remeshing process is important in the long term stability of the
method. In relaxation simulations, once the shape has minimized its
energy there still persists a small velocity field on the membrane.
These velocities, arguably similar to the parasitic velocities
that appear in capillary flows \cite{gmt07,gr07jcp,reusken08,popinet09,asb10}, 
slowly distort the mesh until some sort of instability is
triggered and the simulation diverges. 

\begin{figure}[htp]
\begin{center}
\scalebox{0.8}{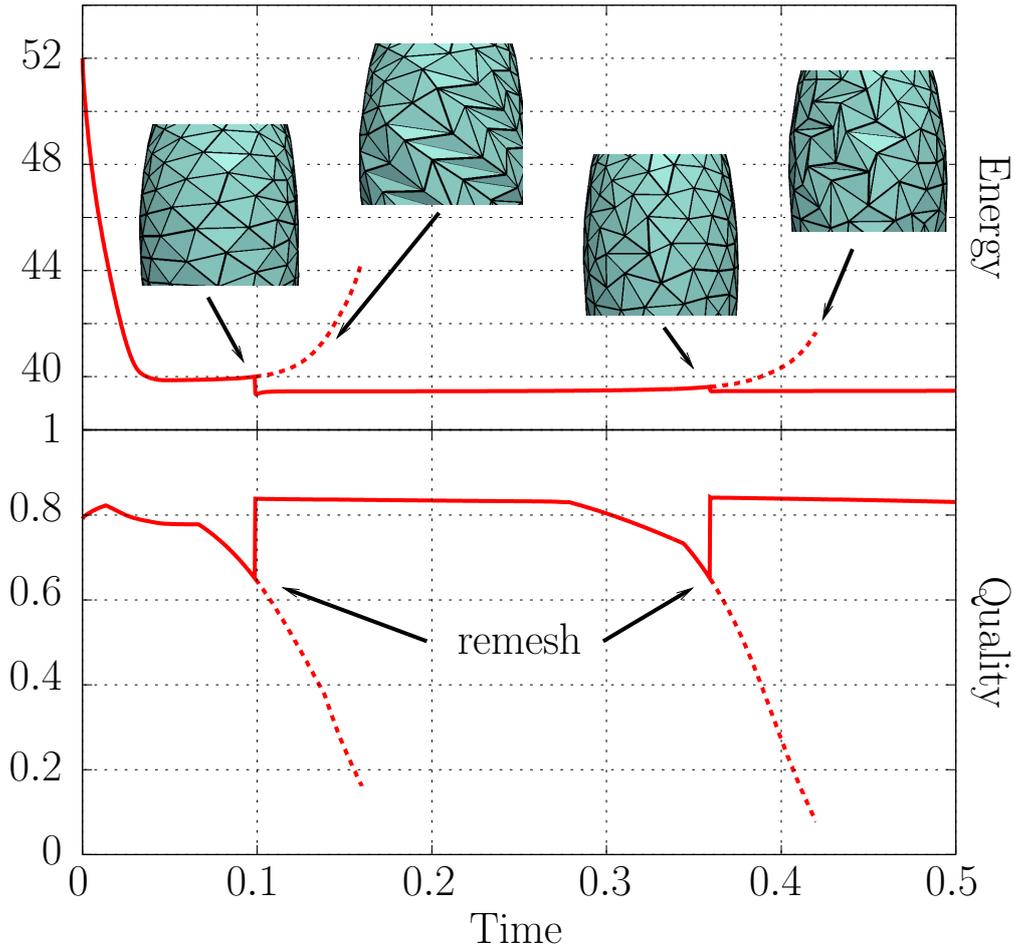}
\end{center}
\caption{Long-term evolution of a relaxing membrane. Plotted are
the elastic energy and mesh quality as functions of time. The state
of rest is not completely achieved and small parasitic velocities
distort the mesh activating the remeshing process. The dotted lines
show the evolution of the variables if remeshing is inhibited,
and the inserts show the unstable distortion pattern.}
 \label{figremesh}
\end{figure}

An attempt to illustrate this phenomenon
is made in Figure \ref{figremesh}. There the evolution of the energy
and of the mesh quality along a relaxation simulation are plotted.
The relaxation should end at $t\simeq 0.06$, with the velocity
going to zero and the membrane remaining forever after in the
equilibrium configuration. One observes, however, that the
quality of the mesh deteriorates steadily and the elastic energy
begins to grow after $t \simeq 0.07$. This behavior, if allowed
to progress, completely pollutes the simulation. The dotted curves
after $t=0.1$ correspond to the evolutions of energy and mesh quality
that would be obtained if the remeshing operation automatically
activated at $t=0.1$ were inhibited. The mesh distortions in these
instabilities are more pronounced in some localized region. The
inserts in Figure \ref{figremesh} show the affected region at
the time of remeshing and sometime later, in a non-remeshed 
simulation. 

After remeshing at $t=0.1$ there is a slight adjustment of the
energy due to the change in mesh and then again a state of
pseudo-rest develops, in which nothing happens other than a
slow distortion of the nodal positions. After time $t\simeq 0.3$
this spurious movement begins to significantly affect mesh
quality and the energy begins to grow again. A new instability
develops quite similar to the one that activated the first
remeshing, leading to a second remeshing at $t\simeq 0.36$.  
The inserts show the critical regions, exhibiting the unstable
distortion pattern.

Remeshing is thus seen to serve not just as a mesh adaptation
strategy, but also as a control mechanism for spurious unstable
distortions.

\subsection{Tethering experiments: Membrane tweezing, dynamical effects}

\subsubsection{Basic description}

\begin{figure}[htp]
\begin{center}
\scalebox{0.8}{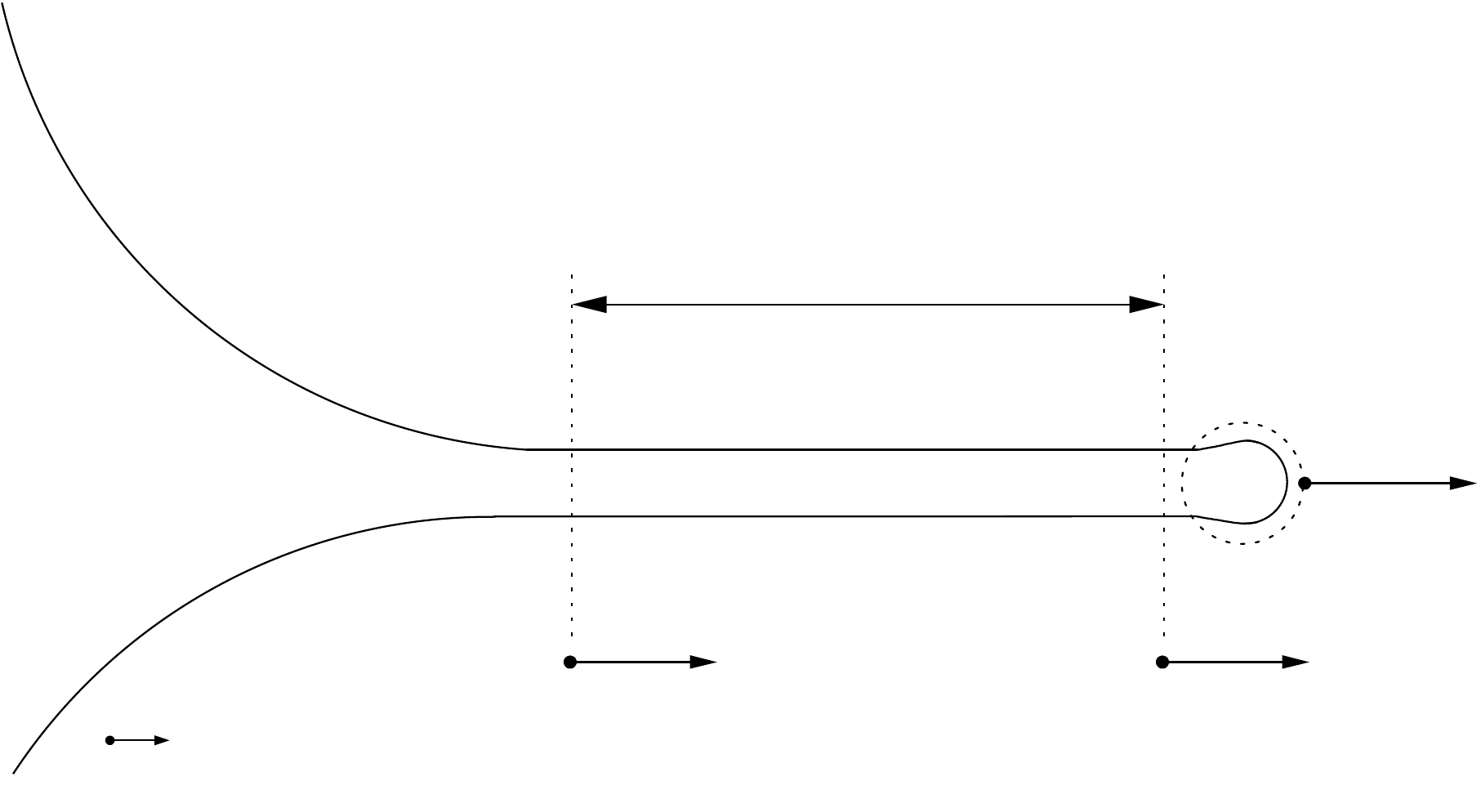}
\end{center}
\caption{Schematics of membrane tethering.}
\label{ffig.tether_scheme}
\end{figure}

A tether develops when a small parcel of the membrane is pulled away.
If a force $F_T$ is applied to the parcel, a structure develops
composed of a head, a cylindrical tube of length $L(t)$
and radius $R(t)$ and the connection to the membrane body
as shown in Figure \ref{ffig.tether_scheme}.

The dynamics of the tether can be understood with the help
of the analytical solution corresponding to a (circular)
cylindrical inextensible membrane. We here go back to
dimensional quantities and consider a membrane of
surface viscosity $\mu$ and Canham-Helfrich's constant 
$c_{\text{\tiny{CH}}}$ that is being pulled from its end by an
external axial force $F_T$. In this particular geometry,
the exact problem admits an analytical solution
with uniform (independent of ${\bf x}$) circumferential
and axial stresses. The exact velocity field is given by
\begin{equation}
\mathbf{u} \,=\,
U_r \, \widecheck{\mathbf{e}}_r
+ \chi\, z \, \widecheck{\mathbf{e}}_z
\end{equation}
with $U_r$
and $\chi$ given by
\begin{eqnarray}
 U_r &=&
-\frac{1}{8 \pi \mu}\left[F_T
\;\;-\;\; 2 \pi c_{\text{\tiny{CH}}}\frac1R
\left(1 + \frac{p\, R^3}{c_{\text{\tiny{CH}}}} \right)\right], \\
\nonumber \\
\chi &=&
\frac{1}{8 \pi \mu R}\left[F_T
\;\;-\;\; 2 \pi c_{\text{\tiny{CH}}}\frac1R
\left(1 + \frac{p\, R^3}{c_{\text{\tiny{CH}}}} \right)\right].
\end{eqnarray}
Neglecting the contribution of the internal pressure $p$,
and noticing that
\begin{equation}
 \frac{dR}{dt} = U_r,
\end{equation}
one arrives at the more tractable equation
\begin{equation} \label{eq.edo_R}
 \frac{dR}{dt} \;=\; \frac{c_{\text{\tiny{CH}}}}{4 \mu}
 \left( \frac{1}{ R_{\text{\tiny{eq}}}}
 - \frac{1}{R(t)} \right).
\end{equation} There exists an {\em equilibrium radius}
$R_{\text{\tiny{eq}}}$ given by
\begin{eqnarray}
 R_{\text{\tiny{eq}}}
 &=& \frac{2 \pi c_{\text{\tiny{CH}}}}{F_T},
\end{eqnarray}
to which the cylinder will tend as $t\,\to\,\infty$. At equilibrium,
the surface pressure $\pi_{\text{s}}$ takes the value
$$
\pi_{\text{s,eq}}=
-\,\frac{F_T}{4\pi R_{\text{\tiny{eq}}}}=
-\,\frac{F_T^2}{8\pi^2 c_{\text{\tiny{CH}}}}.
$$
Further, the final decay when $R \simeq R_{\text{\tiny{eq}}}$
must have the asymptotic behavior
\begin{equation}
R(t)=R_{\text{\tiny{eq}}}+~C~\exp\left ( -\frac{c_{\text{\tiny{CH}}}}{4\,\mu\,R_{\text{\tiny{eq}}}^2}
~t \right ).
\end{equation}
The characteristic relaxation time is 
$$
\mathcal{T} = \frac{4\,\mu\,R_{\text{\tiny{eq}}}^2}{c_{\text{\tiny{CH}}}}
=\frac{16\,\pi^2\,\mu\,c_{\text{\tiny{CH}}}}{F_T^2}.
$$
For $t$ much greater than $\tau$ the tether is expected to be
at equilibrium following a rigid-body translation along the
line of $F_T$. The material deforms to take the shape of a cylinder
in the region to the left of point ``$b$'' in Figure \ref{ffig.tether_scheme},
which is approximately fixed in space (the ``beginning'' of the tether).
Once the material enters the tether it simply moves at constant
velocity along it. The ``end'' of the tether (point ``$e$'') moves
at a constant velocity $U_T$ determined by a balance between the
applied force and the viscous stresses at the connection region between
the tether and the membrane body.

Going back to non-dimensional quantities, the equilibrium radius
and the tether relaxation time are given by
\begin{equation}
R_{\text{\tiny{eq}}}=\frac{2\,\pi}{F_T}
\qquad
\mbox{and}
\qquad
\mathcal{T}=\frac{16\,\pi^2}{F_T^2}.
\label{eqreqtau}
\end{equation}

\subsubsection{Numerical tweezers}

We have implemented numerical tweezers as a model for the external
surface force ${\bf f}$. Each numerical tweezer has a radius
$r_T$, which is fixed in time, while the position of its center
follows a path described by the vectors ${\bf x}_T^{0}$,
${\bf x}_T^{1}$, etc. 

Given a point ${\bf x}$ in $\mathbb{R}^3$, the tweezer's
penetration at point ${\bf x}$ and time $t_{n+1}$, 
denoted by $w^{n+1}({\bf x})$, is defined as 
$$
w^{n+1}({\bf x})~\stackrel{\mbox{\tiny{def}}}{=}~r_T-\|{\bf d}_T^{n+1}({\bf x})\|,
$$
where
$$
{\bf d}_T^{n+1}({\bf x})~\stackrel{\mbox{\tiny{def}}}{=}~ {\bf x}-{\bf x}_T^{n+1}.
$$
The repulsive force that the tweezer exerts on ${\bf x}$ depends
exponentially on $w^{n+1}({\bf x})$, according to
$$
{\bf f}^{n+1}({\bf x})~=~k_T\,\frac{e^{w^{n+1}({\bf x})/\ell_T}}{\|{\bf d}_T^{n+1}({\bf x})\|}
~{\bf d}_T^{n+1}({\bf x}).
$$
In an exact setting, this force would be integrated over 
${\bf x}\,\in\,\Gamma^{n+1}$.
Unfortunately, this force is needed at the time of computing
${\bf u}_h^{n+1}$ through (\ref{eqvaruh}), and thus the integral 
is performed over $\Gamma^n$. One could replace ${\bf f}^{n+1}$ by ${\bf f}^n$
in (\ref{eqvaruh}), but the following approximation has much more
stable behavior:
\begin{equation}
{\bf f}^{n+1}({\bf x})~=~
~k_T\,\frac{e^{w^{n+1}({\bf x})/\ell_T}}{\|{\bf d}_T^{n}({\bf x})\|}
~{\bf d}_T^{n}({\bf x})~
-~k_T\,\frac{\delta t\,e^{w^{n}({\bf x})/\ell_T}}{\ell_T\,\|{\bf d}_T^{n}({\bf x})\|^2}
\,\left ({\bf d}_T^{n}({\bf x})\otimes {\bf d}_T^{n}({\bf x})\right )
~{\bf u}_h^{n+1}({\bf x}).
\label{eq57}
\end{equation}
Notice that the last term in (\ref{eq57}) is an implicit linearization
that must be added to the matrix arising from the left-hand side
of (\ref{eqvaruh}). The two parameters $k_T$ and $\ell_T$ are given the
values $10^5$ and $r_T/50$, respectively.

The numerical tweezer can be moved specifying either the velocity or
the total force exerted on the membrane. 
In the former case, the update rule is simply
$$
{\bf x}_T^{n+1}={\bf x}_T^{n}+ \delta t\,{\bf U}_T^{n+1},
$$
where ${\bf U}_T$ is the specified tweezer velocity.
In cases where the force is specified, a simple proportional feedback controler
was implemented that adjusts ${\bf U}_T^{n+1}$ so as to keep the force
at the target value.

\begin{figure}[htp]
\begin{center}
\scalebox{0.6}{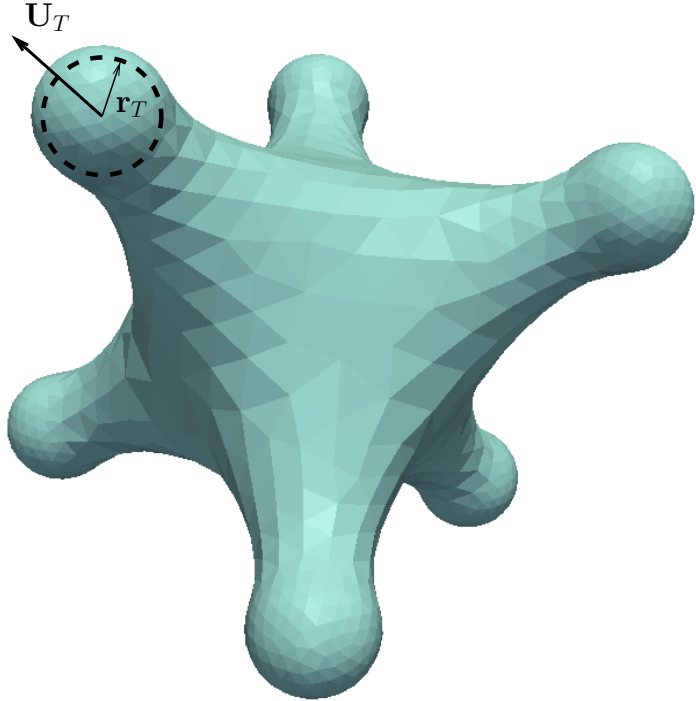}
\end{center}
\caption{An illustration of the tweezing strategy. In this case six
tweezers are simultaneously acting on a membrane.}
 \label{figtweezer}
\end{figure}

Figure \ref{figtweezer} illustrates our tweezing strategy,
depicting a situation in which six tweezers
are simultaneously pushing a membrane outwards, from within.

\subsubsection{Dynamical effects in tweezing}

The viscous-flow model presented here can deal, unlike gradient-flow models,
with dynamical effects in excursions away from equilibrium. An illustrative
application is the analysis of velocity effects in tether development
carried below.

Let us consider an equilibrium oblate shape corresponding to a reduced
volume $v=0.6666$. A tweezer of constants
$$
r_T=0.04, \qquad k_T=10^5,\qquad \ell_T=\frac{r_T}{50}
$$
is placed at ${\bf x}_T(t=0)$ on the interior of the surface, close to it, and is 
moved at constant velocity $U_T$ along the outward normal direction.
The simulations are run until a time $T$ such that the
tweezer displacement $D_T=\|{\bf x}_T(T)-{\bf x}_T(0)\|=U_T\,T = 1.0$,
meaning a non-dimensional displacement of the tweezer of 1.0 for
all cases.
We have observed that it is necessary to reduce the time step for
the tweezing simulations. Specifically, $\delta t$ is now chosen as
\begin{equation}
\delta t~=~\frac{1}{4}~\delta t^*(h_{\min})
\label{eqdttweezer}
\end{equation}
and the results are confirmed by re-running the simulation with 
one half of this value.

Since the goal is to consider just the interaction of one tweezer,
the rigid motions are filtered out by Lagrange multipliers just as
in the free relaxation cases.

\begin{figure}[htp]
\begin{center}
\scalebox{0.6}{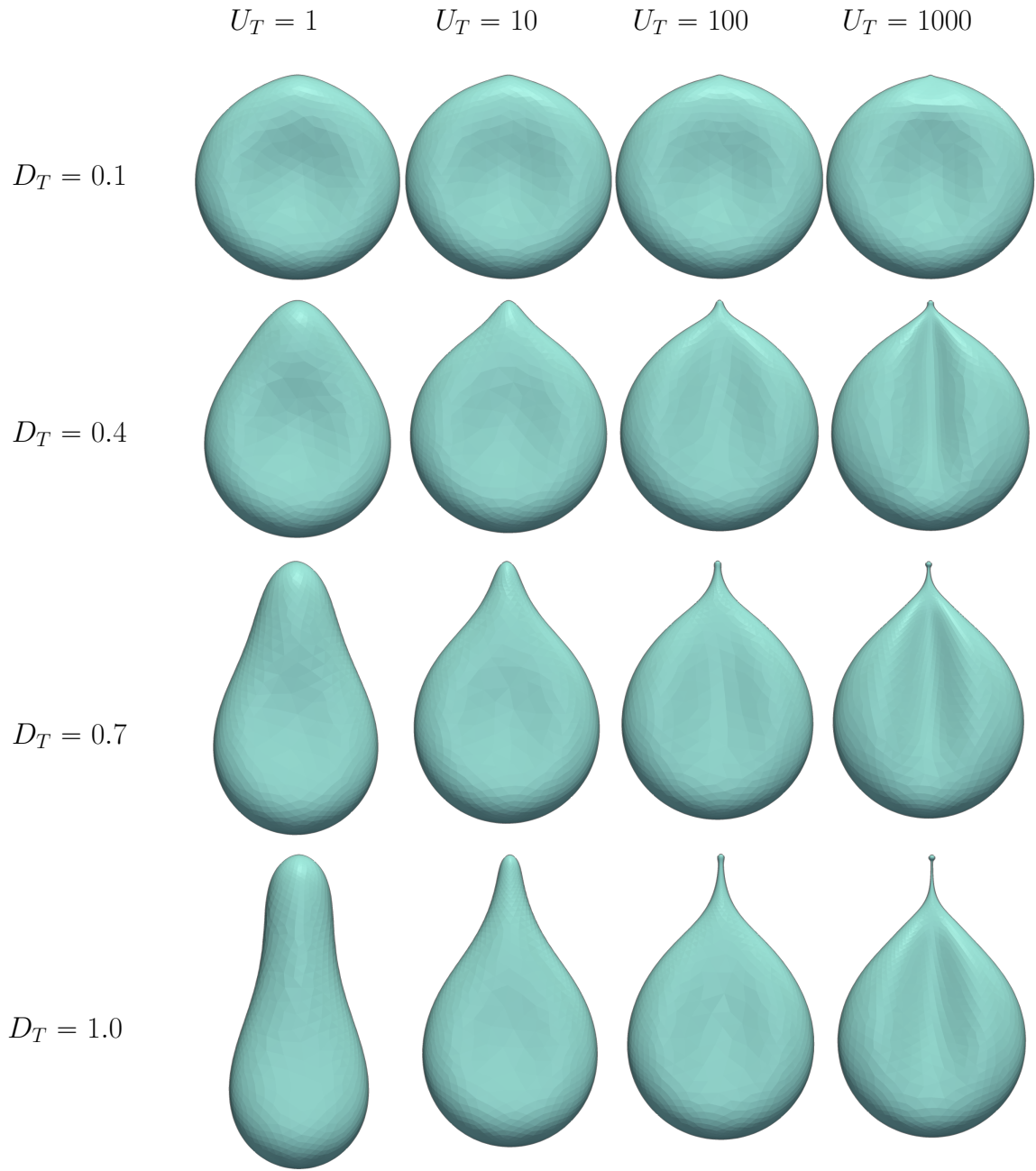}
\end{center}
\caption{Dynamical effects in tweezing. Deformation of an equilibrium
oblate shape with $v=0.6666$ by a tweezer of radius $r_T=0.04$ moving
outwards at constant velocity $U_T$ (vertically). Shown
are the membrane shapes for four values of the tweezer's displacement
$D_T=0.1$, 0.4, 0.7 and 1.0, and for $U_T=1$, 10, 100 and 1000.
}
\label{figdynamictweezing}
\end{figure}

The resulting membrane shapes at different positions of the tweezer,
indicated by its displacement $D_T$, and for several values of the
tweezer velocity $U_T$ are shown in Figure \ref{figdynamictweezing}.

The leftmost column of the figure corresponds to the smallest
velocity, $U_T=1$. In this case, the membrane deforms almost
quasistatically, without showing any localized response
at the tweezer's location. As the velocity is increased to $U_T=10$
one begins to ``see'' the tweezer pushing outwards from within
the membrane. But it is only for $U_T=100$ and $U_T=1000$
that the small size of the tweezer ($r_T=0.04$) becomes apparent
and the tweezer produces a tethering-like deformation. 

Considering just the bottom row of Figure \ref{figdynamictweezing},
for which the tweezer position is exactly the same (and the
center of mass of the membrane too, thanks to the rigid-motion
filtering), the ability of the proposed method to capture
velocity-dependent deformations of the membrane is evident.

\subsubsection{Tether dynamical equilibrium}

We now assess the ability of the method to correctly predict the
dynamical equilibrium of the tether. For that purpose, we take a
tether formed by applying a tweezer force of $F_T=400$ and
suddenly change the force. Two runs were performed, in one of
them $F_T$ is changed to $500$ and in the other to $600$.
The time is redefined to be zero at the time of the force change.
The mesh is adapted and remeshed using $c_h=1/2$, which is rather
coarse ($h\,\simeq\,R_{\text{\tiny{eq}}}/2$). The time step is adjusted
according to (\ref{eqdttweezer}).

\begin{figure}[htp]
\begin{center}
\scalebox{0.88}{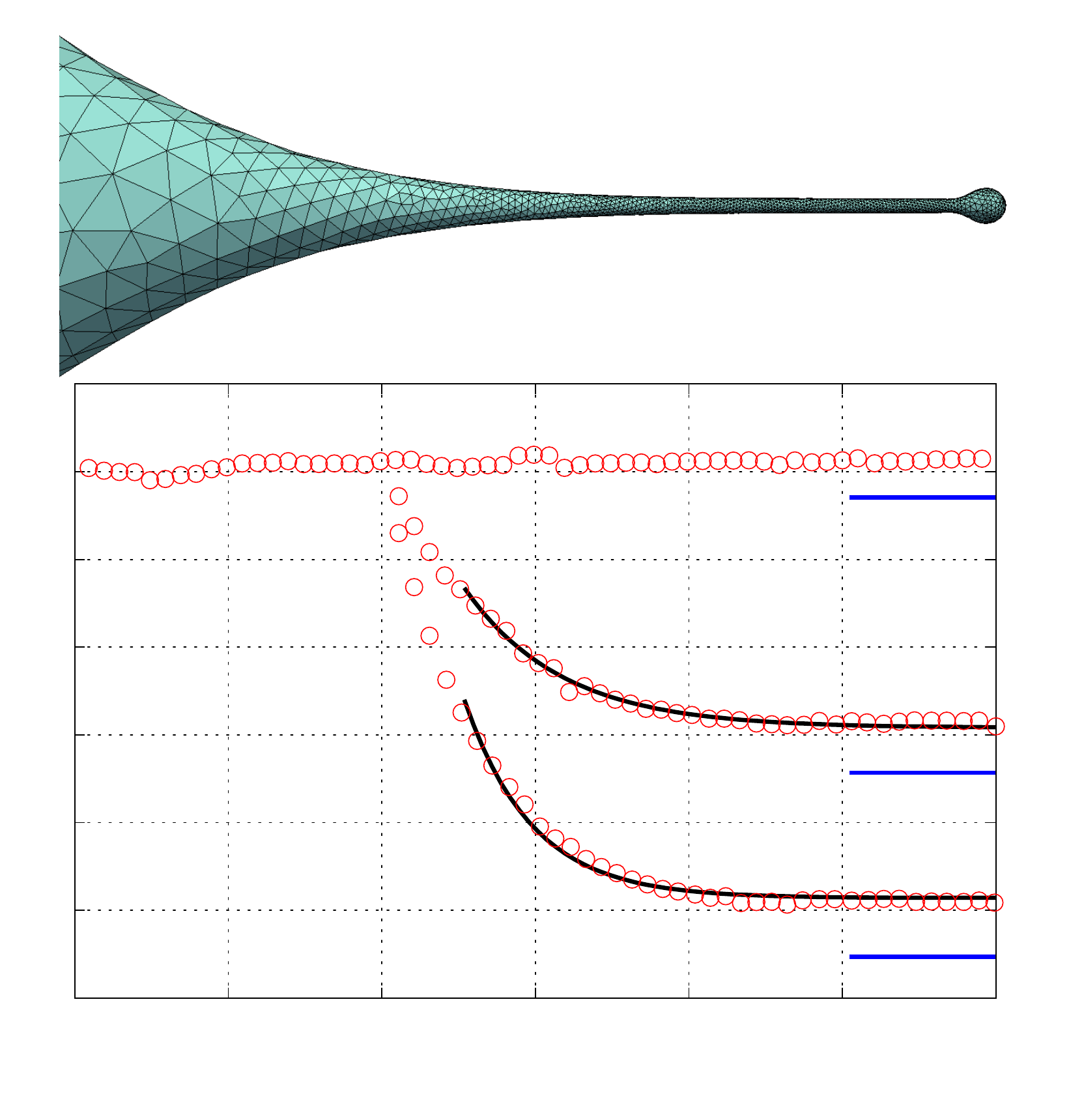}
\end{center}
\caption{Response of the tether's radius to a sudden change in the
tweezer's force. From a tether in dynamical equilibrum at $F_T=400$,
the force is changed to $F_T=500$ or $F_T=600$ at $t=0$. The circles correspond 
to the numerical results obtained with the proposed method. 
The black lines correspond to an exponential adjustment with
characteristic time given by (\ref{eqreqtau}). On the right 
the exact equilibrium radius (as given by (\ref{eqreqtau}) is
shown with a short blue segment.
}
 \label{figrequilibrium}
\end{figure}

By post-processing the mesh it is possible to compute the
radius of the tether as a function of time, as shown in Figure
\ref{figrequilibrium}. The initial exact equilibrium radius
is $R_{\text{\tiny{eq}}}(F_T=400)=0.0157$, which is reasonably
approximated by the method despite the mesh being quite coarse.

After changing the force to $F_T=500$, the tether's radius shrinks
to a value of approximately 0.013, which is a good approximation to
$R_{\text{\tiny{eq}}}(F_T=500)=0.0126$. Further, the evolution towards
the new radius is in good agreement with an exponential of the
form $a\,e^{-t/\mathcal{T}}+b$, where $\mathcal{T}=6.32\times 10^{-4}$ is given by (\ref{eqreqtau}),
as shown by a continuous line in the figure.

A similar procedure is conducted for the change to $F_T=600$,
for which the exact values are $R_{\text{\tiny{eq}}}(F_T=600)=0.0105$
and $\mathcal{T}=4.39\times 10^{-4}$.

The relaxation towards the equilibrium radius is seen to agree
quite well with the analytical solution (though better for $F_T=500$ than
for $F_T=600$), and the equilibrium
radius itself is predicted with an error of about 5\%. 
This error level is reasonable, considering that there are just
about twelve elements in the tether's circumference.

\subsubsection{Complex tweezing}

This last section reports on a more complex tweezing experiment
which aims at testing the robustness of the proposed method.
Starting from a spherical membrane of area $4\,\pi$ (i.e.,
taking as $R_0$ the radius of the initial sphere), six
independent tweezers of radius $r_T=0.1$ 
move radially outwards with $U_T=100$
acting upon it. The tweezers initial positions are the
intersection of the membrane with the six cartesian semi-axes.

Though the imposed values of area, $\mathcal{A}^*$, and
volume, $\mathcal{V}^*$, were introduced as constants, in this
experiments they are set as specified functions of time,
\begin{eqnarray}
\mathcal{A}^*(t)&=& 4\pi\,+\,400\,t,\\
\mathcal{V}^*(t)&=& \frac{4\pi}{3}\,-\,50\,t. 
\end{eqnarray}
A time-dependent enclosed volume may result from a variable
osmotic pressure in the fluid that surrounds the membrane,
while a time-dependent area may result from the incorporation
of lipids to the membrane.

Along the simulation, the time step was continually adjusted
according to
$$
\delta t~=~0.105\,h_{\min}^2
$$
and the remeshing procedure was applied automatically,
with $c_h=0.5$.

\begin{figure}[htp]
\begin{center}
\scalebox{0.75}{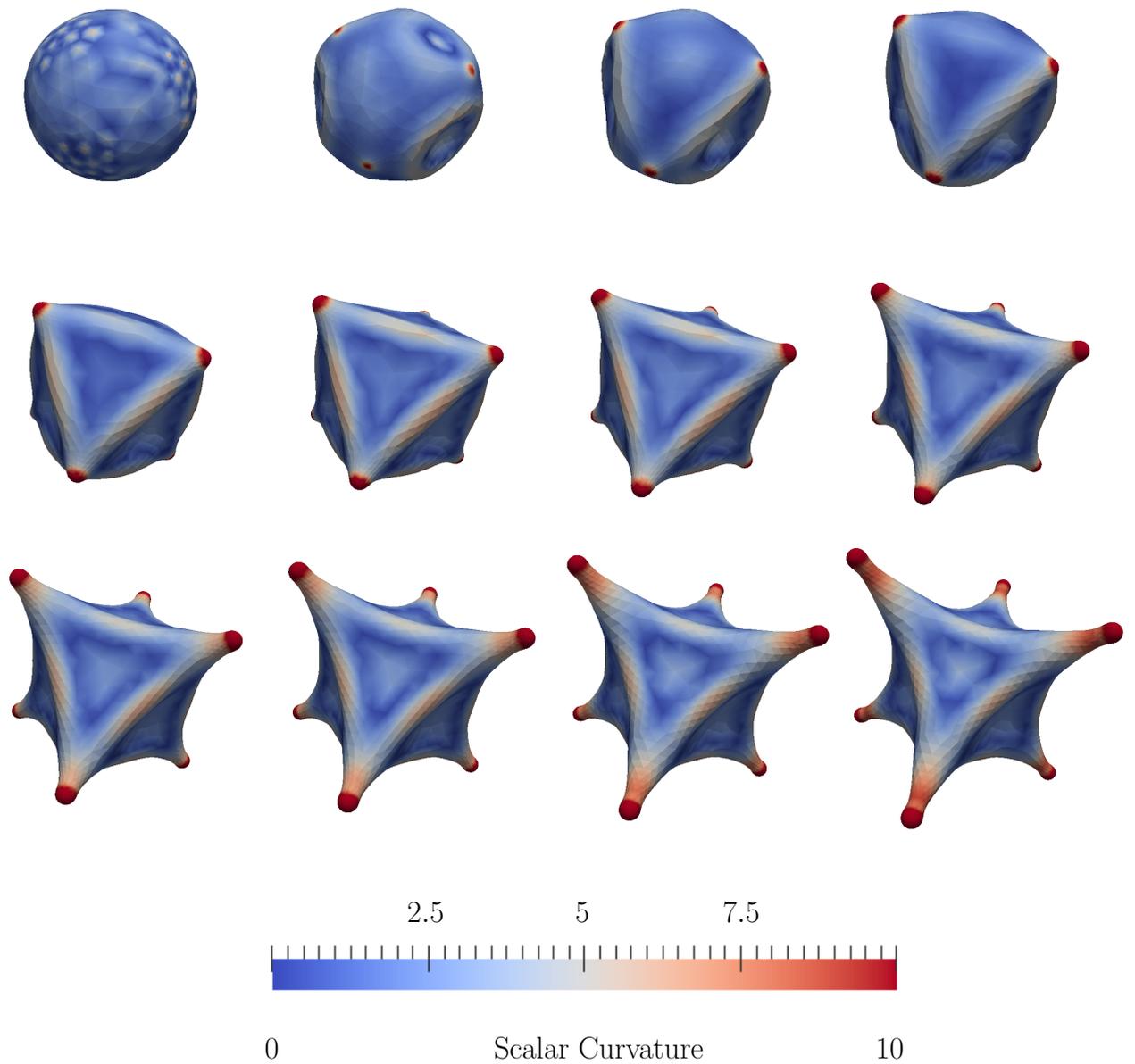}
\end{center}
\caption{Evolution of the membrane's shape along the
``complex tweezing'' simulation. Shown are snapshots of
the shape at equispaced time intervals of $0.8\times 10^{-3}$
time units. The shapes are shaded according to the value
of the scalar curvature.}
 \label{figcomplex}
\end{figure}

A picture of the membrane's evolution can be seen in Figure
\ref{figcomplex}. The tweezers are seen to ``emerge'' from
the sphere first deforming the membrane into an approximate octahedron
(at time $\sim 0.005$) and then further stretching the octahedron
into a star-like shape. Though there exist mechanisms that may
create protrusions such as those in Figure \ref{figcomplex} in
actual cells or lipid vesicles \cite{fml97,sars13}, this case does
not attempt to model a specific physical phenomenon.

\begin{figure}[htp]
\begin{center}
\scalebox{0.8}{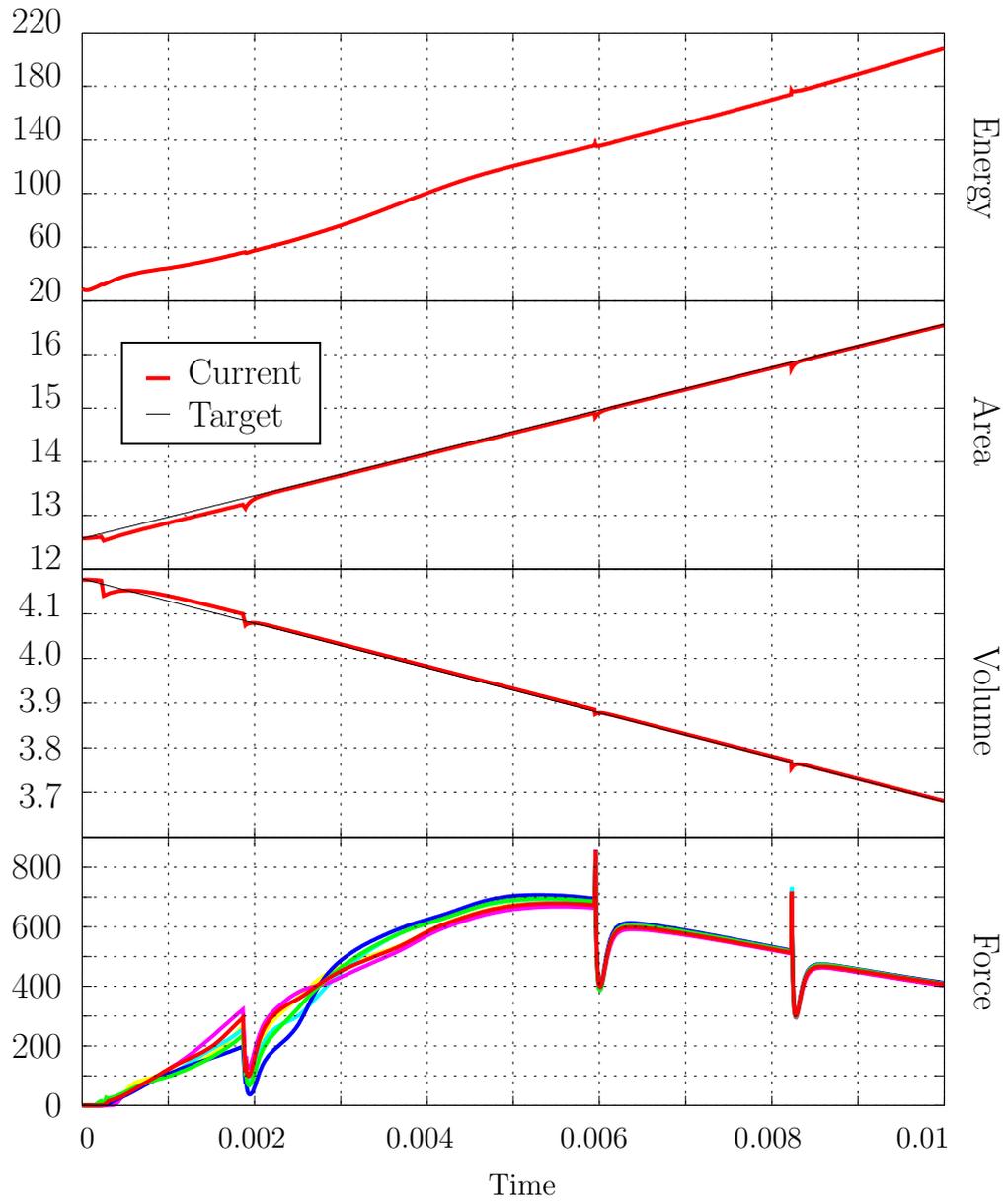}
\end{center}
\caption{Plots of energy, area, volume and tweezer forces
(the six of them) as functions of time as obtained in the
``complex tweezing'' simulation.}
 \label{figcomplex2}
\end{figure}

In Figure \ref{figcomplex2} plots of several variables of the
simulation can be found. The energy is seen to increase
monotonically along the deformation, with the area and volume
following their target values $\mathcal{A}^*(t)$ and
$\mathcal{V}^*(t)$ quite closely. The forces exerted by each
of the six tweezers are also plotted in Figure \ref{figcomplex2}.
They differ from one another until at $t\simeq 0.005$ the
membrane tightens and all the tweezers start behaving alike.

Notice the strong perturbations introduced by remeshing, which
are the result of slight changes in the penetration of each
tweezer by the interpolatory construction of the new mesh.
The algorithm is able to recover itself from these strong
perturbations quite rapidly.

\begin{figure}[htp]
\begin{center}
\scalebox{0.8}{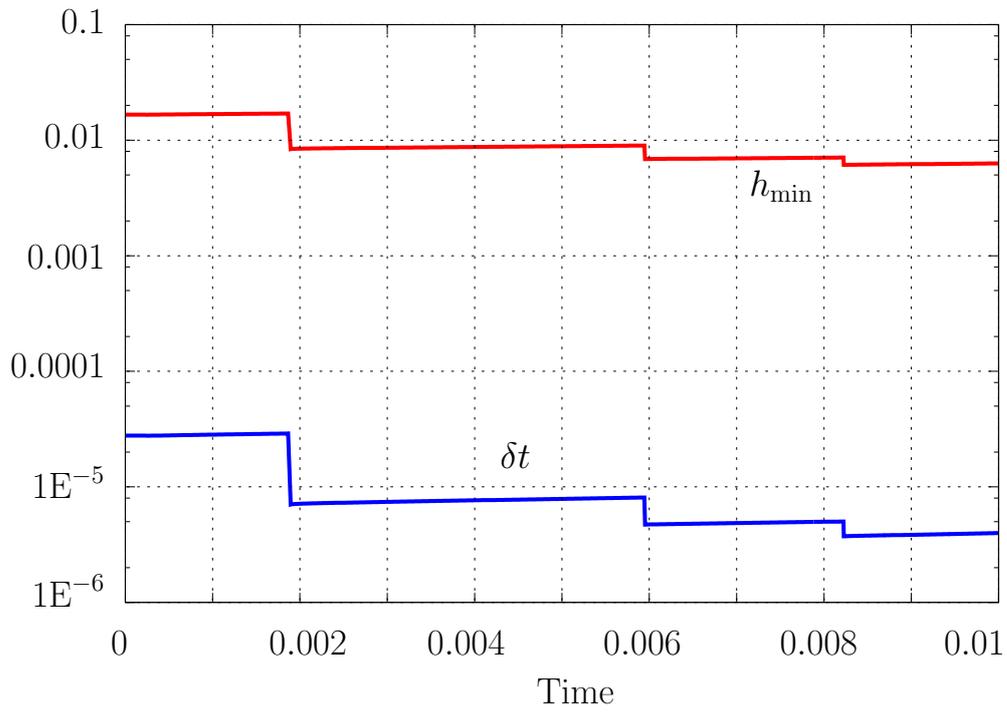}
\end{center}
\caption{Plots of $\delta t$ and $h_{\min}$ as functions of
time in the ``complex tweezing'' simulation.}
 \label{figcomplex3}
\end{figure}

\begin{figure}[htp]
\begin{center}
\scalebox{0.8}{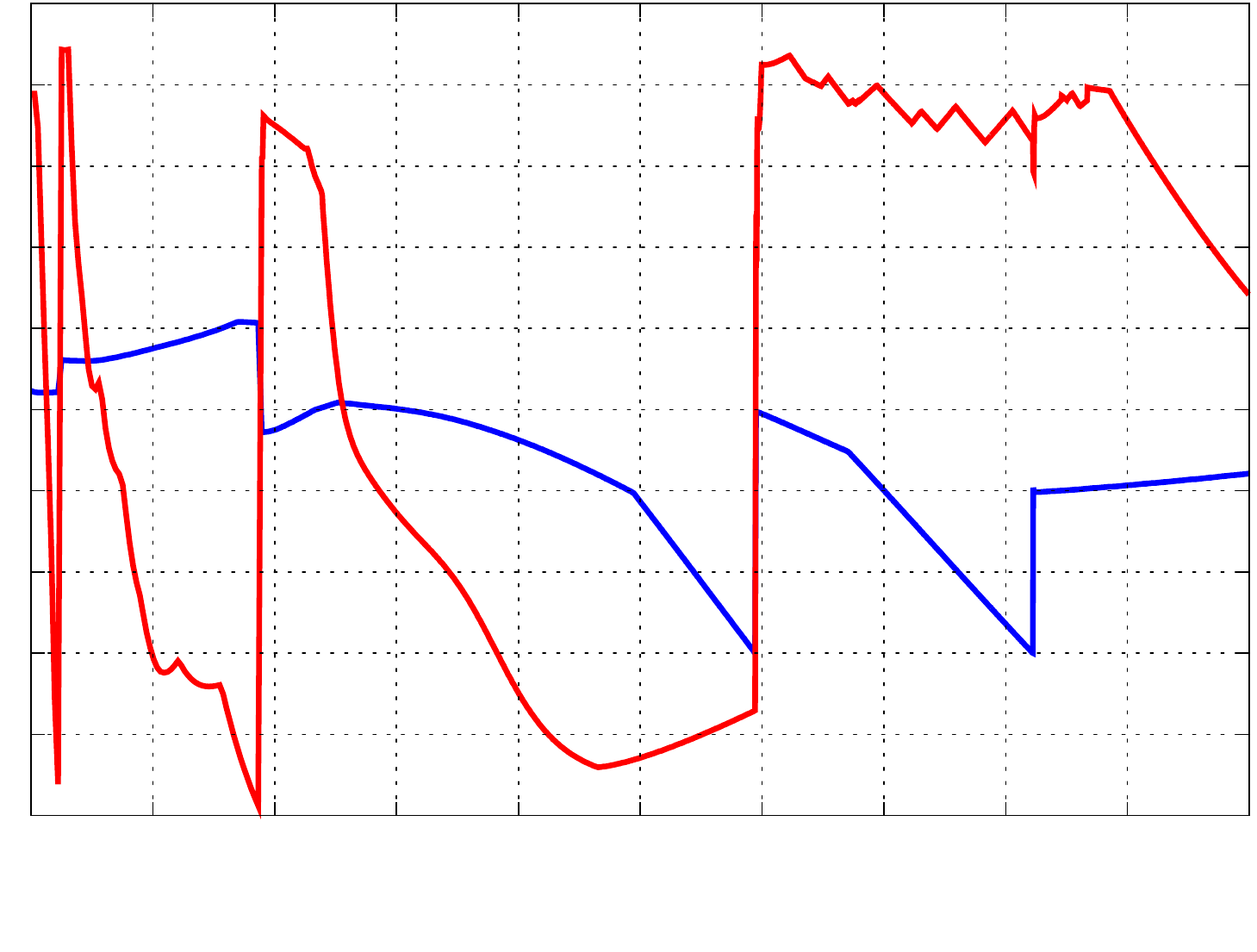}
\end{center}
\caption{Plots of $Q_{\mbox{\tiny{shape}}}$ and $Q_{\mbox{\tiny{size}}}$
as functions of time in the ``complex tweezing'' simulation.}
 \label{figcomplex4}
\end{figure}

Finally, let us provide some computational data of this
simulation. The time step and the minimum element size
$h_{\min}$ are plotted as functions of time in Figure
\ref{figcomplex3}. Also shown are the shape and size qualities of the
mesh, $Q_{\mbox{\tiny{shape}}}$ and $Q_{\mbox{\tiny{size}}}$, as functions
of time in Figure
\ref{figcomplex4}. The initial mesh consists of 2160 elements and 1082
nodes, while the final one consists of 5128 elements and
2566 nodes. The complete simulation comprises 1525 time steps,
which take 29 minutes on an i7-based laptop at 2.8 GHz.
The linear system is solved by LU factorization using
the MUMPS package \cite{MUMPS:1,MUMPS:2}, with a memory requirement of 2 GBytes.

\section{CONCLUSIONS}

In this contribution, we have introduced a fully discrete semi-implicit
finite element scheme for the simulation of viscous membranes with
bending elasticity of the Canham-Helfrich type. 
The membrane is discretized
by a surface mesh made up of planar triangles, over
which a mixed formulation (velocity-curvature) is
built with $P_1$ interpolants for all fields.
Two stabilization terms are incorporated in the
discrete formulation: The first one stabilizes
the inextensibility constraint by a pressure-gradient-projection
scheme \cite{cb97}, the second couples curvature and velocity to
improve temporal stability \cite{bansch01}.
The volume constraint is handled by a Lagrange multiplier
(which turns out to be the internal pressure), and an
analogous strategy is used to filter out rigid-body motions.
Feedback controllers are used to avoid drifting from imposed
values of enclosed volume and total area.
The nodal positions are updated in a Lagrangian manner
and automatic remeshing strategy maintains suitable refinement
and mesh quality throughout the simulation.

The method has been numerically assessed through extensive relaxation
and tweezing experiments. For the latter, a specific virtual
tweezing algorithm was devised. It has been shown
that the proposed method is convergent and robust, though
with a severe (of order $h^2$) stability restriction on the time 
step for which a practical estimate was derived. This stability restriction
is the main difficulty in the applications of the algorithm, since
it makes thousands of time steps necessary for
the simulation of relatively simple membrane motions.

Another difficulty still encountered, though currently avoided by
quality-based automatic remeshing, is the existence of small persistent
velocities at the numerical equilibrium which slowly but continually
deteriorate the mesh quality.

\section*{ACKNOWLEDGMENTS}

The authors gratefully acknowledge the financial support received 
from S\~ao Paulo Research Foundation (FAPESP, grants. no. 2011/01800-5,
2012/14481-8 and 2012/23383-0) and from the Brazilian National
Research and Technology Council (CNPq).

\bibliographystyle{plain}
\bibliography{biblio}

\end{document}